\documentclass[a4paper,11pt]{article}

\usepackage{mathrsfs} 

\usepackage{amsthm}
\usepackage{hyperref}
\usepackage{amssymb}
\usepackage{latexsym}
\usepackage{amsmath}
\usepackage{amsfonts}
\usepackage{mathabx}
\usepackage[dvips]{graphicx}
\usepackage[utf8]{inputenc}
\usepackage{comment}
\usepackage[LGR,T1]{fontenc}%
\usepackage[french,english]{babel}
\usepackage{url}
\usepackage{enumerate}
\usepackage{hyperref}
\usepackage{color}
\usepackage{subfig,graphicx}

\usepackage{fancyhdr}
\newtheorem{Theorem}{Theorem}[part]
\newtheorem{Definition}{Definition}[part]
\newtheorem{Proposition}{Proposition}[part]

\newtheorem{Lemma}{Lemma}[part]
\newtheorem{Corollary}{Corollary}[part]
\newtheorem{Remark}{Remark}[part]
\newtheorem{Example}{Example}[part]

\makeatletter \@addtoreset{equation}{section}

\@addtoreset{Definition}{section}

\@addtoreset{Theorem}{section}

\@addtoreset{Proposition}{section}

\@addtoreset{Property}{section}

\@addtoreset{Assumption}{section}

\@addtoreset{Corollary}{section}

\@addtoreset{Lemma}{section}

\@addtoreset{Remark}{section}

\@addtoreset{Example}{section}

\newenvironment{bluetext}{\color{blue}}{\ignorespacesafterend}

\newcommand{\cA}{\mathcal{A}}

\newcommand{\cD}{\mathcal{D}}
\newcommand{\cE}{\mathcal{E}}
\newcommand{\cF}{\mathcal{F}}

\newcommand{\cK}{\mathcal{K}}
\newcommand{\cL}{\mathcal{L}}
\newcommand{\cM}{\mathcal{M}}
\newcommand{\cN}{\mathcal{N}}

\newcommand{\cP}{\mathcal{P}}

\newcommand{\cS}{\mathcal{S}}
\newcommand{\cT}{\mathcal{T}}

\newcommand{\cV}{\mathcal{V}}

\newcommand{\cY}{\mathcal{Y}}
\newcommand{\cZ}{\mathcal{Z}}

\renewcommand{\P}{\mathbb{P}}

\newcommand{\R}{\mathbb{R}}

\def \proof{{\noindent \bf Proof. }}
\def \eproof{\hbox{ }\hfill$\Box$}

\newcommand{\ud}{\mathrm{d}}

\newcommand{\1}{{\bf 1}}
    
\newcommand{\set}[1]
    {\ensuremath{\{ #1 \}}}
    
\newcommand{\HP}[1] 
    {\ensuremath{\mathscr{H}^{#1}}}

\newcommand{\esp}[1]{\ensuremath{\mathbb{E} \!\! \left[#1\right] }}

\renewcommand{\Xi}[1]{X_{i #1}}

\title{Numerical approximation of singular Forward-Backward SDEs}

\author{
  Jean-Fran{\c{c}}ois CHASSAGNEUX\thanks{UFR de Math{\'e}matiques \& LPSM, Universit{\'e} de Paris, B{\^a}timent Sophie Germain, 8 place Aur{\'e}lie Nemours, 75013 Paris, France ({\tt chassagneux@lpsm.paris})}, Mohan YANG\thanks{UFR de Math{\'e}matiques \& LPSM, Universit{\'e} de Paris, B{\^a}timent Sophie Germain, 8 place Aur{\'e}lie Nemours, 75013 Paris, France ({\tt myang@lpsm.paris})}}

\begin{document}
\maketitle

\begin{abstract}
In this work, we study the numerical approximation of a class of singular fully coupled forward backward stochastic differential equations. These equations have a  degenerate forward component and non-smooth terminal condition. They are used, for example, in the modeling of carbon market \cite{carmona2013singular} and are linked to scalar conservation law perturbed by a diffusion. Classical FBSDEs methods fail to capture the correct entropy solution to the associated quasi-linear PDE. We introduce a splitting approach that circumvent this difficulty by treating differently the numerical approximation of the diffusion part and the non-linear transport part. Under the structural condition guaranteeing the well-posedness of the singular FBSDEs \cite{carmona2013singularconservation}, we show that the splitting method is convergent with a rate $1/2$. We implement the splitting scheme combining non-linear regression based on deep neural networks and conservative finite difference schemes. The numerical tests show very good results in possibly high dimensional framework.
\end{abstract}
\vspace{5mm}


\vspace{5mm}

\noindent {\bf MSC Classification (2020):}  65C30, 65C20, 65M12, 60H35.

\section{Introduction}
In this work, we study the approximation of a class of singular fully coupled Forward Backward Stochastic Differential Equations (FBSDE). Let $(\Omega,\cF,\P)$ be a stochastic basis supporting a $d$-dimensional Brownian motion $W$ and $T>0$ a terminal time. We denote by $(\cF_t)_{t \ge 0}$ the filtration generated by the Brownian motion (augmented and completed). The singular FBSDE system, with solution $(P_t,E_t,Y_t,Z_t)_{0\le t\le T}$, has the following form:
\begin{align} \label{eq singular fbsde}
\left\{
\begin{array}{rcl}
dP_t & = & b(P_t)dt + \sigma(P_t) \ud W_{t}\\
dE_t &=  &\mu(Y_t,P_{t})\ud t\\
dY_t &=  & Z_t \cdot \ud W_t
\end{array}
\right.
\end{align}
The function $b:\R^d \rightarrow \R^d$, $\sigma:\R^d \rightarrow \cM_d$\footnote{\textcolor{black}{To alleviate the notation, we assume that $P$ and $W$ have the same dimension and the coefficient functions of $P$ are time-homogeneous. Note however that $\sigma$ will not be assumed to be uniformly elliptic, which allows to consider a dimension  of $P$ as time and to embed the case of different dimension for $P$ and $W$ in our framework.}}, where $\cM_d$ is the set of $d \times d$ matrices on $\R$, and $\mu : \R\times\R^d \rightarrow \R$ are Lipschitz-continous.
\noindent These equations have been introduced in \cite{carmona2013singular} as models for carbon emission market. They can model, more generally, cap-and-trade scheme used by government to limit the emission of certain pollutant. In these models, $Y$ is the price of a pollution permit, $E$ is the cumulative emission of the pollutant and $P$ represents some state variables influencing the emission (demand, energy prices etc.). The coefficient $\mu$ is naturally decreasing in the $y$-variable. The terminal condition is given by $\phi(E_T,P_T)$, where $\phi : \R \times \R^d \rightarrow \R$ is a measurable function, non-decreasing in its $E$-variable and Lipschitz continuous in the 
$P$-variable. In its simplest form, it is given typically by:
\begin{align}\label{eq term condition}
e \mapsto \phi(e) = \1_{\set{e > \Lambda}}\,, \; \Lambda > 0.
\end{align} 
The constant $\Lambda$ appears as a cap on emissions set by the regulator. The shape given in \eqref{eq term condition} translates the fact that a penalty (here set to one) is paid if the emission are above the regulatory cap at $T$.\\
We observe that \eqref{eq singular fbsde} has a forward one dimensional $E$-component of bounded variation and a backward component with an irregular terminal condition \eqref{eq term condition}. This renders the mathematical analysis of the FBSDE system difficult. Nevertheless, the well-posedness and main features of \eqref{eq singular fbsde} have been thoroughly studied in \cite{carmona2013singularconservation}, see also Section \ref{subse wellposedness} below. Notably, the authors of \cite{carmona2013singularconservation} prove existence and uniqueness of the solution to \eqref{eq singular fbsde} but show at the same time that the terminal condition can only be attained in the following weak sense:
\begin{align}\label{eq weak terminal condition}
\1_{(\Lambda,+\infty)}(E_T) \le Y_T \le \1_{[\Lambda,+\infty)}(E_T) \;,
\end{align}
using to simplify the presentation at this point the terminal function \eqref{eq term condition}.
Their study is based on the celebrated markovian representation of $Y$ as 
\begin{align} \label{eq de decoupling field}
Y_t = \cV(t,P_t,E_t), \;\text{ for }\; t<T,
\end{align} 
and the careful analysis of the property of $\cV$, where $\cV$, known as the \emph{decoupling fields}, is solution to a quasilinear PDE. As mentioned in \cite{carmona2013singularconservation}, the FBSDE system can be seen as a random perturbation of a scalar conservation law. The behavior at the terminal time is reminiscent of shocks appearing in conservation law. Let us note that the markovian representation breaks down at $T$ as indicated by \eqref{eq weak terminal condition}. Moreover, the function $\cV$ is only locally Lipschitz-continuous on $[0,T)$:
\begin{align}\label{eq gradient explosion}
|\cV(t,p,e)-\cV(t,p',e')| \le c_1|p-p'| + \frac1{c_2(T-t)}|e-e'|\,,
\end{align}
for some constants $c_1,c_2>0$, $(p,p',e,e') \in \R^d\times\R^d\times\R\times\R$.

\begin{figure} 
\centering
\subfloat[$\sigma=0.01$]{\label{fig a}\includegraphics[width=.3\linewidth]{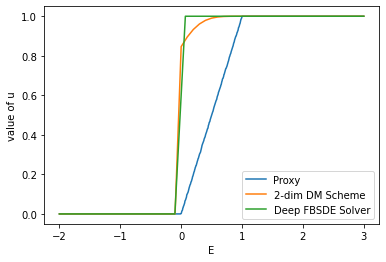}}\hfill
\subfloat[$\sigma=0.3$]{\label{fig b}\includegraphics[width=.3\linewidth]{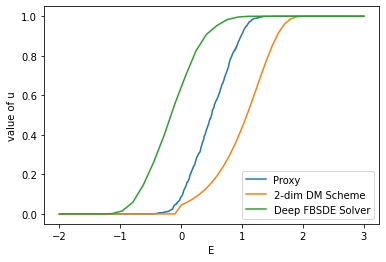}}\hfill 
\subfloat[$\sigma=1.0$]{\label{fig c}\includegraphics[width=.3\linewidth]{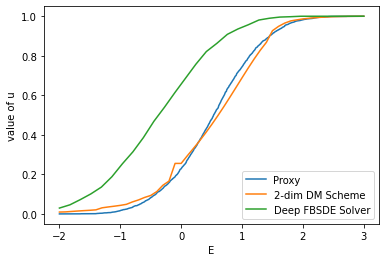}}
\caption{Comparison of $e \mapsto \cV(0,0,e)$ obtained by Deep FBSDE Solver and Delarue-Menozzi Scheme (DM Scheme) to the proxy, for different level of volatility. The methods fail to reproduce correctly the proxy.}\label{fig main figure}
\end{figure}

The application to carbon market is also a key motivation for our numerical study here: efficient numerical simulation of the price $Y$ would allow to calibrate properly  the model to market data and validate its efficiency in practice.  A first approach for the numerical approximation of $Y$ or $\cV$ would be to use PDE methods, and this is suggested in \cite{howison2012risk}. However, in the economic applications we have in mind, the dimensionality of the process $P$ prevents generally the use of these methods. In order to work on problems in moderate dimension, say 5 to 10, some probabilistic methods could be introduced. 
Probabilistic schemes have already been designed for FBSDEs  and one could be tempted (as we were) to use the already known methods to tackle the numerical approximation of \eqref{eq singular fbsde}. 
In \cite{bender2008time}, the authors use a Picard Iteration method to decouple the FBSDE system and then obtain an approximation of $\cV$ by performing iteratively linear regression. Unfortunately, this  method has only been shown to be convergent in the case of Lipschitz coefficient and for small coupling between the forward and backward part (or equivalently small time horizon), see \cite{bender2008time} for details. Recently, machine learning methods have been considered for BSDEs approximation, especially for their applicability in very high-dimensional setting. In particular, \cite{han2020convergence} has analysed the \emph{deep BSDE solver} introduced in \cite{Han8505} again in the setting of small coupling. In
\cite{delarue2006forward}, a grid algorithm  is introduced where the decoupling is obtained by a predictor: there, the time horizon or the coupling is arbitrary but the diffusion coefficient of the forward process must be uniformly elliptic.  As observed, the FBSDE system under study is degenerated in the $E$-component and the terminal condition is not Lipschitz, so that none of the known methods for FBSDEs are proved to be convergent in the setting of \eqref{eq singular fbsde}.
Moreover, the above methods fail, in practice, to approximate correctly the solution to \eqref{eq singular fbsde}. 
\noindent To empirically illustrate this fact, we consider the following toy model borrowed from \cite{carmona2013singular}:
\begin{Example}[Linear model]\label{ex lin model}
\label{ex linear emission}
\begin{align}
dP_t& = \sigma\ud W_{t} \label{eq toy model P}\\
dE_t &= \left(\frac1{\sqrt{d}}\sum_{\ell=1}^dP^\ell_t-Y_t\right)\ud t \label{eq toy model E}\\
dY_t &=  Z_t \cdot \ud W_t \label{eq toy model Y}
\end{align}
with terminal function given by \eqref{eq term condition} and 
where $W$ is a $d$-dimensional Brownian motion and $\sigma > 0$.
\end{Example}
By using a change of variable,  this $d+1$ dimensional model can be reduced to a one dimension model. Indeed, from  \cite[Proposition 6]{carmona2013singular}, there exists $\nu \in C^{1,2}([0,T),\R)$, solution to
\begin{align}\label{eq pde change variable}
\partial_t \nu - \nu \partial_\xi \nu + \frac{\sigma^2(T-t)^2}{2} \partial^2_{\xi\xi} \nu = 0 \; \text{ and }\; \nu(T,\xi) = \phi(\xi)\,.
\end{align}
By essentially applying Ito's formula (see the proof of  \cite[Proposition 7]{carmona2013singular} for details on proving \eqref{eq pde change variable}), one obtains that $\cV(t,p,e) = \nu(t, e + (T-t)\frac{1}{\sqrt{d}}\sum_{\ell=1}^d p_\ell)$. 
This observation allows us to use efficient methods to solve \eqref{eq pde change variable} numerically and to compare them to numerical solutions obtained by ``classical'' FBSDE scheme. 
In particular, we use a probabilistic method  studied in \cite{bossy1997stochastic} using interacting particle system, a class of mean field SDE, to obtain a ``proxy'' for $e \mapsto \cV(0,0,e)$, see also \cite{bossy1996convergence,jourdain2002probabilistic}.
\\
Going back to the approximation of \eqref{eq toy model P}-\eqref{eq toy model E}-\eqref{eq toy model Y} by ``classical'' FBSDEs methods, we first note that, in \cite[Chapter 4]{fbsdespringerbrief}, the authors report an application of the Bender-Zhang scheme \cite{bender2008time}. The main issue is then that the Picard iteration does not converge to a single limit.
Next, we have tested the Delarue-Menozzi scheme \cite{delarue2006forward} 
and the \emph{deep FBSDE solver} \cite{han2020convergence} for different value of $\sigma$, the results are given in Figure \ref{fig main figure}.

Except maybe for the Delarue-Menozzi scheme in Figure 
\ref{fig c}, the methods fail clearly to approximate the correct solution $\cV(0,0,\cdot)$. The problem comes from the nonlinear transport part of the equation in this degenerate setting.  Indeed, the methods seem unable to recover the correct weak entropy solution. This is particularly clear on Figure \ref{fig a}, where the level of noise is extremely small and the correct solution is almost the solution to the inviscid Burger's equation. This leads us to introduce a new method to approximate the FBSDE system \eqref{eq singular fbsde}. 

As already mentioned, in the socio-economic applications, the dimension of the $P$-variable is generally large. On the contrary, the $E$-variable is constrained to be of dimension one. We note also that approximating the dynamics of the $P$-variable corresponds to approximating simply a diffusion process, which can be easily done. To take into account these key differences in the two variables, we follow a splitting approach to compute numerically the solution $\cV$. On a discrete time grid, we iterate backward in time, a diffusion operator where the $E$-variable is fixed to capture the effect of the $P$-dynamics in \eqref{eq singular fbsde}, and a transport operator where the $P$ variable is fixed to capture the effect of $E$-dynamics in \eqref{eq singular fbsde}. Our main theoretical result, see Theorem \ref{th conv res theo scheme}, proves that this scheme is convergent at a rate $\frac12$ with respect to the time step. Our analysis is done under the minimal assumption used in \cite{carmona2013singularconservation} to obtain existence and uniqueness of the solution $\cV$. One of the main difficulty encountered is therefore due to the gradient explosion at the end of the time interval \eqref{eq gradient explosion}.

Then, we propose various implementations of the splitting scheme. They have however a common structure: given a discrete transport operator, the diffusion part is computed by means of probabilistic methods. The overall scheme is a then a sequence of (non linear) regression in the high dimensional space where lives the approximation of $\cV$ with respect to the $E$-variable. 
In our numerical experiments, we consider approximations of the transport operator by conservative finite difference methods (Lax-Friedrichs scheme or Upwind scheme), see e.g. \cite{leveque1992numerical}. As we do not always have access to a proxy for the tested models, we introduce an alternative implementation of the splitting scheme to validate our numerical results. It combines  a particle approximation of the transport operator with a tree based regression for the diffusion operator. We validate empirically both approaches on Example \ref{ex lin model} for which we have a one dimensional proxy at hand. We then test models with no equivalent one-dimensional PDE but whose $d+1$ dimensional version can be reduced to $2$-dimensional specification, see Section \ref{se num exp} for details. The tree-based algorithm is then used as a proxy as it is very efficient in low-dimension (we test dimension 4) for the models under consideration. When combining feedforward neural networks to compute the regression step and finite difference scheme for the transport step, we show that our splitting procedure can compute precisely and in reasonable amount of time the solutions of $10+1$ dimensional models. 


\vspace{2mm}
The rest of the paper is organised as follows. In Section \ref{se splitting}, we first recall key properties of the theoretical solution. We then introduce the splitting approach and prove the convergence of the splitting scheme.
In Section \ref{se num scheme}, we present a regression method for the splitting scheme at a theoretical level, which uses a grid approximation of the transport operator. We then introduce various implementation of the transport operator and a neural network approximation for the regression part. We finally present various numerical experiments to validate the efficiency of our method in practice.
%

\subsubsection*{Notation.}

In the  following we will  use the following spaces  
\begin{itemize}
\item 
For  fixed $0\le a<b<+\infty$ and $I=[a,b]$ or $I=[a,b)$, { $\mathcal{S}^{2,k}(I)$ is the
set of  $\mathop{}\!\mathbb{R}^k$-valued c\`adl\`ag\footnote{French acronym for right continuous with left limits.} $\mathcal{F}_t$-adapted processes $Y$, s.t.\\[-4mm]  
\begin{align*}
\|Y\|_{\mathcal{S}^2}^2:= \mathbb{E}\left[{\sup_{t \in I} |Y_t|^2 } %
\right] <\infty .
\end{align*}
Note that we may omit the dimension and the terminal date in the norm notation as this will be clear from the context. $\mathcal{S}^{2,k}_\mathrm{c}(I)$ is the subspace of process with continuous sample paths.} 

\item For fixed $0\le a<b<+\infty$, and $I=[a,b]$, we denote by { $\mathcal{H}^{2,k}(I)$ the set
of $\mathop{}\!\mathbb{R}^k$-valued progressively measurable processes $Z$, such that \\[-4mm] 
\begin{equation*}
\|Z\|_{\mathcal{H}^2}^2 := \mathbb{E} \left[{\int_I |Z_t|^2 dt} \right] 
<\infty.  
\end{equation*}%
}
\end{itemize}


%

\noindent For $\varphi:\R^d \times \R \rightarrow \R$, measurable and non-decreasing in its second variable,
the functions $\varphi_{-}$ and $\varphi_{+}$ are the left and right continuous versions, respectively defined, for $(p,e) \in \mathbb{R}^{d} \times \mathbb{R}$, by,
	\begin{align}
	\begin{aligned}
	\label{phi minus plus definition}
	\varphi_{-}(p,e) &= \sup_{e^{\prime} < e} \varphi(p, e^{\prime}) \\
	\varphi_{+}(p,e) &= \inf_{e^{\prime} > e} \varphi(p, e^{\prime}).
	\end{aligned}
	\end{align}
Moreover, we denote by $\| \cdot \|_\infty$ the essential supremum:
\begin{align*}
\| \varphi \|_\infty = \mathrm{esssup}_{(p,e)\in \R^d \times \R} |\varphi(p,e)| \;.
\end{align*}

%
%



\section{A splitting scheme}
\label{se splitting}

In this section, we introduce  a theoretical splitting scheme to compute the solution of the singular FBSDEs. This scheme consists into iterating a ``diffusion step'' and a ``transport step'' on a discrete time grid 
$$\pi := \set{ 0=: t_0 <\dots<t_n<\dots<t_N := T},$$
where $N$ is a positive integer. For latter use, we denote by $|\pi| := \max_{0\le n < N}(t_{n+1}-t_n)$. 

Before defining the splitting scheme for the system \eqref{eq singular fbsde}, we recall some key theoretical properties of the solution obtained in \cite{carmona2013singularconservation}, with slight extensions for the case of P-dependent terminal condition in \cite{jean2020modelling}.  The rest of the section is then dedicated to the proof of an upper bound for the convergence rate of the splitting scheme in terms of $|\pi|$. This is our main theoretical result, given in Theorem \ref{th conv res theo scheme}. Numerical implementations are presented in the next section.
 
\subsection{Well-posedness and properties of singular FBSDEs}
\label{subse wellposedness}

We first introduce two classes of functions, that will be useful in the sequel.  
The terminal condition function for \eqref{eq singular fbsde} will belong to the first one. 
\begin{Definition}\label{de cK}
Let $\cK$ be the class of functions $\phi:\R^d\times\R \rightarrow [0,1]$ such that $\phi$ is $L_\phi$-Lipschitz  in the first variable for some $L_\phi>0$ and non-decreasing in its second variable, namely
\begin{align}
|\phi(p,e)-\phi(p',e)| &\le L_\phi|p-p'| \quad\text{ for all }\quad (p,p',e) \in \R^d\times\R^d\times\R\;,
\\
\phi(p,e') &\ge \phi(p,e) \quad\text{ if }\quad e' \ge e\;,
\end{align} 
and moreover  satisfying,
\textcolor{black}{
\begin{align}\label{eq condition on phi}
\sup_e \phi(p,e) = 1  \;\text{ and }\; \inf_e \phi(p,e) = 0 \quad\text{ for all }\quad p \in \R^d\;.
\end{align}
}
\end{Definition}
\noindent \textcolor{black}{Note that the bounds given in  \eqref{eq condition on phi} are motivated by our main application, but up to a rescalling they can be arbitrary changed.}
We now introduce a class of admissible coefficient functions, for which the singular BSDE is well-posed, see Theorem \ref{th existence uniqueness one-period} below. This class will be also useful to define the splitting scheme.
\\
We consider three positive constants $L, \ell_1$ and $\ell_2$.

\begin{Definition} \label{de class coef} Let $\cA$ be the class of functions 
$B:\R^d \rightarrow \R^d$, $\Sigma:\R^d \rightarrow \cM_d$, $F:\R\times \R^d\rightarrow \R$ which are  $L$-Lipschitz continuous functions. Moreover, \textcolor{black}{$F$ is strictly decreasing in $y$} and satisfies, for all $p \in \R^d$,
\begin{align}\label{eq conservation law}
\ell_1 |y-y'|^2 \le (y-y')(F(y',p)-F(y,p)) \le \ell_2 |y-y'|^2.
\end{align}
\end{Definition}

\noindent \textbf{Standing assumptions:} From now on, we assume that $(b,\sigma,\mu) \in \cA$, recalling \eqref{eq singular fbsde}.%

\begin{Theorem}[Proposition 2.10 in \cite{carmona2013singularconservation}, Proposition 3.2 in \cite{jean2020modelling}] \label{th existence uniqueness one-period} 
Let $\tau > 0$, $(B,\Sigma,F) \in \cA$ and $\Phi \in \cK$.

\noindent Given any initial condition $(t_0,p,e) \in [0,\tau) \times \mathbb{R}^{d} \times \mathbb{R}$, there exists a unique progressively measurable 4-tuple of processes $(P^{t_0,p,e}_t, E^{t_0,p,e}_t, Y^{t_0,p,e}_t, Z^{t_0,p,e}_t)_{t_0 \leq t \leq \tau} \in \mathcal{S}^{2,d}_{\mathrm{c}}([t_0,\tau]) \times \mathcal{S}^{2,1}_{\mathrm{c}}([t_0,\tau]) \times \mathcal{S}^{2,1}_{\mathrm{c}}([t_0,\tau)) \times \mathcal{H}^{2,d}([t_0,\tau])$ satisfying the dynamics
\begin{align}  \label{general single period fbsde with p}
\begin{aligned}
\ud P^{t_0,p,e}_t &= B(P^{t_0,p,e}_t) \ud t + \Sigma(P^{t_0,p,e}_t) \ud W_t, & P^{t_0,p,e}_{t_0} &= p \in \mathbb{R}^{d}, \\
\ud E^{t_0,p,e}_t &= F(P^{t_0,p,e}_t, Y^{t_0,p,e}_t) \ud  t, & E^{t_0,p,e}_{t_0}
&= e \in \mathbb{R}, \\ 
\ud Y^{t_0,p,e}_t &=  Z^{t_0,p,e}_{t} \cdot  \ud {W_t}, & &
\end{aligned}
\end{align}
and such that
	\begin{align}
	\label{relaxed terminal condition with p}
	\mathbb{P} \left[ \Phi_{-}(P^{t_0,p,e}_\tau,E^{t_0,p,e}_\tau) \leq \lim_{t \uparrow \tau}Y^{t_0,p,e}_t \leq \Phi_{+}(P^{t_0,p,e}_\tau, E^{t_0,p,e}_\tau) \right] = 1.
	\end{align}
The unique \emph{decoupling field} defined by
\begin{align*}
[0,\tau) \times \mathbb{R}^{d} \times \mathbb{R} \ni (t_0,p,e) \rightarrow w(t_0,p,e) = Y^{t_0,p,e}_{t_0} \in \R
\end{align*}	
is continuous and satisfies
\begin{enumerate}
		\item{\label{v lipschitz constant e}For any $t \in [0,\tau)$, the function $w(t,\cdot, \cdot)$ is $1/(l_1(\tau-t))$-Lipschitz continuous with respect to $e$,}
		\item{\label{v lipschitz constant p}For any $t \in [0,\tau)$, the function $w(t,\cdot, \cdot)$ is $C$-Lipschitz continuous with respect to $p$, where $C$ is a constant depending on $L$, $\tau$ and $L_{\phi}$ only.}
		\item{\label{v condition terminal} Given $(p,e) \in \R^d \times \R$, for any family $(p_t, e_t)_{0 \leq t < \tau}$ converging to $(p,e)$ as $t \uparrow \tau$, we have
		\begin{align}\label{eq v condition terminal}
		\Phi_-(p,e) \le \liminf_{t \rightarrow \tau}w(t,p_t,e_t) \le \limsup_{t \rightarrow \tau}w(t,p_t,e_t) \le \Phi_+(p,e)\,.
		\end{align}
		\item{\label{v belongs to K}}
		\textcolor{black}{For any $t \in [0,\tau)$, the function $w(t,\cdot, \cdot) \in \cK$.}
}
\end{enumerate}
\end{Theorem}

\noindent Using the previous result, we define the following operator associated to \eqref{eq singular fbsde}.
\begin{Definition}\label{de ope true solution}
We define the operator $\Theta$ by
\begin{align}
(0,{\infty})\times \cK \ni (h,\psi) \mapsto \Theta_h(\psi) = v(0,\cdot) \in \cK\,
\end{align}
where $\upsilon$ is the decoupling field given in Theorem \ref{th existence uniqueness one-period} with parameters $\tau = h$, $B=b$, $\Sigma = \sigma$, $F=\mu$ and $\Phi = \psi$.
\end{Definition}
We also deduce from Theorem \ref{th existence uniqueness one-period}  that $(\Theta_t)_{0<t}$ is a semi-group of non-linear operators. 
In particular, we observe that $\cV(0,\cdot) := \Theta_T(\phi) = \prod_{0 \le n  < N}\Theta_{t_{n+1}-t_n}(\phi)$, recall \eqref{eq de decoupling field}.

\vspace{2mm}
\noindent The following result arises from the proof of the previous Theorem, see \cite{jean2020modelling}.
\begin{Corollary}[Approximation result]\label{co smooth approx of v}
Let $\tau > 0$, $(B,\Sigma,F) \in \cA$ and $\Phi \in \cK$. Let $(\phi^k)_{k \ge 0}$ be a sequence of smooth functions belonging to $\cK$ and converging pointwise towards $\phi$ as k goes to $+\infty$. For $\epsilon>0$, consider then $w^{\epsilon,k}$ the solution to:
\begin{align}
\label{value function pde}
\partial_t u +F(p,u) \partial_{e} u + \cL_p u + \frac12\epsilon^2 (\partial^2_{ee}u+ \Delta_{pp}u) = 0 \; \text{ and } \; u(\tau,\cdot) = \phi^k
\end{align}
where $\Delta_{pp}$ is the Laplacian with respect to $p$, and $\cL_p$ is the operator
\begin{align}
\cL_p(\varphi)(t,p,e) = \partial_p\varphi (t,p,e)B(p) + \frac{1}{2} \mathrm{Tr} \left[ A(p) \partial^{2}_{pp} \right](\varphi)(t,p,e),
\end{align}
with \textcolor{black}{$\partial_p$ denotes the Jacobian with respect to $p$},  and  $A = \Sigma \Sigma^{\top}$, where $\top$ is the transpose and $\partial^{2}_{p p}$ is the matrix of second derivative operators. (For later use, we define
$ 
\cL^\epsilon := \cL_p  + \frac12\epsilon^2 (\partial^2_{ee}+ \Delta_{pp})\;.
$
)

Then the functions $w^{\epsilon,k}$ are $C^{1,2}$ (continuously differentiable in $t$ and twice continuously differentiable in both $p$ and $e$) and $\lim_{k \rightarrow \infty} \lim_{\epsilon \rightarrow 0} w^{\epsilon,k} = w$ where the convergence is locally uniform in $[0,\tau)\times\R^d \times \R$.
Moreover, for all $k,\epsilon$, $w^{k,\epsilon}(t,\cdot) \in \cK$.
\end{Corollary}


\subsection{Scheme Definition} 

Let us first introduce the transport step where the diffusion part is frozen.
\begin{Definition}[Transport step] \label{de transport step}
 We set
 $$(0,\infty) \times \cK \ni (h,\psi) \mapsto \cT_h(\psi) = \tilde{v}(0,\cdot) \in \cK$$ 
 where $\tilde{v}$ is the decoupling field defined in Theorem \ref{th existence uniqueness one-period} with parameters
 $\tau = h$, $B =0$, $\Sigma = 0$, $F = \mu$ and terminal condition $\Phi = \psi$.
\end{Definition}

\noindent In the definition above, $\tilde{v}(\cdot)$ is the unique entropy solution to 
\begin{align} \label{eq de cons law associated with transport}
\partial_t w + \partial_e ( \mathfrak{M}(p,w) ) = 0\,, \quad \text{ where } \mathfrak{M}(p,y)=\int^y_0 \mu(p,\upsilon) \ud \upsilon\;,\;
0\le y \le 1,
\end{align}
and $\tilde{v}(h,\cdot)=\psi$. 
We will use this fact in the numerical section.


\vspace{2mm} We now introduce the diffusion step, where conversely, the $E$ - process is frozen to its initial value.
\begin{Definition}[Diffusion step]\label{de diffusion step}
 We set
 $$(0,\infty) \times \cK \ni (h,\psi) \mapsto \cD_h(\psi) = \bar{v}(0,\cdot) \in \cK$$ 
 where $\bar{v}(0,\cdot)$ is the decoupling in Theorem \ref{th existence uniqueness one-period} with parameters
 $\tau = h$, $B =b$, $\Sigma = \sigma$ , $F = 0$ and terminal condition $\Phi = \psi$.
\end{Definition}
\noindent Observe that, for $t \in [0,h)$,
\begin{align}\label{eq other de bar v}
\bar{v}(t,p,e) = \esp{\psi(P_h^{t,p},e)}\; \text{ and }  \bar{v}(t,\cdot) \in \cK\;.
\end{align}

\vspace{2mm}
\noindent We can now define the theoretical scheme on $\pi$ by a backward induction.
\begin{Definition}[Theoretical splitting scheme] \label{de th splitting theo}
We set 
\begin{align*}
(0,\infty)\times \cK \ni (h,\psi) \mapsto \cS_h(\psi) := \cT_h \circ \cD_h(\psi) \in \cK. 
\end{align*}

\noindent For $n \le N$, we denote by $u^\pi_n$ the solution of the following backward induction on $\pi$:
\\
- for $n=N$, set $u^\pi_N := \phi$,
\\
- for $n < N$, $u^\pi_n = \cS_{t_{n+1}-t_n}(u^\pi_{n+1})$.
\end{Definition}
\noindent The $(u^\pi_n)_{0 \le n \le N}$ stands for the approximation of the decoupling field $\cV(t,\cdot)$ for $t \in \pi$. Moreover, we observe, from the property of $\cT$ and $\cD$, that 
\begin{align}\label{eq u in K}
u^\pi_n \in \cK, \text{ for all } 0 \le n \le N.
\end{align}

\subsection{Convergence analysis}


Our main theoretical result concerning the splitting is the following
\begin{Theorem}\label{th conv res theo scheme}
Under our \emph{standing assumptions}, the following holds 
\begin{align*}
\int_\R|\cV(0,p,e)-u^\pi_0(p,e)|\ud e \le C T(1+|p|^2) \sqrt{|\pi|}\,,
\end{align*}
for a positive constant $C$.
\end{Theorem}

\noindent The proof of the Theorem is postponed to the end of the section. It is classically based  on the study of the scheme's stability and its truncation error.

\subsubsection{Truncation error}

We need to compare, for $\psi \in \cK$, $\Theta_h(\psi)$ and $\cS_h(\psi)$, $h>0$, to assess the truncation error. As already mentioned, the true solution $\cV$ has minimal locally Lipschitz regularity and it exhibits a gradient explosion in the $E$-variable near the terminal time $T$. In the proof below, we thus need to consider smoothed version of the  decoupling fields introduced in the  definition of $\Theta$, $\cT$,  $\cD$ and $\cS$.
\\
First of all, for a given $\psi \in \cK$, we consider a smooth approximation sequence $\psi^k$ as in Corollary \ref{co smooth approx of v}. In particular, $v^{k,\epsilon}$ is the smooth approximation of the decoupling field $v= \Theta_h(\psi)$  in Definition \ref{de ope true solution} and the associated FBSDEs, for $0\le t \le h$, $Y^{k,\epsilon}_t = v^{k,\epsilon}(t,E^{k,\epsilon}_t,P^\epsilon_t)$
\begin{align}
P^\epsilon_t &= p + \int_0^t b(P^\epsilon_s)\ud s + \int_0^t \sigma(P^\epsilon_s)\ud W_s + \epsilon W'_t \;,
\\
E^{k,\epsilon} &= e +\int_0^t \mu(Y^{k,\epsilon}_s,P^\epsilon_s) \ud s + \epsilon B_t\;,
\end{align}
where $(W',B)$ is a Brownian motion independent from $W$. Note that for the reader's convenience, we omit the dependence upon the starting point $(0,p,e)$ in the FBSDEs notation.
The convergence of $v^{k,\epsilon}$ to $v$ is given in Corollary \ref{co smooth approx of v}.
\\
We also need to consider a smooth version of $\cS_h(\psi)$, that we define now: 
\begin{enumerate}
\item  for $0 \le t \le h$,  set:
\begin{align} \label{eq de bar v}
\bar{v}^{k,\epsilon}(t,p,e) = \esp{\psi^k(P^\epsilon_{h-t},e)}
\end{align}
\item then, $\tilde{v}^{k,\epsilon}$ is the decoupling of the following FBSDE, for all $p \in \R^d$, $e \in \R$: 
\begin{align} 
\ud \tilde{Y}^{k,\epsilon}_t &= \tilde{Z}^{k,\epsilon}_t \ud B_t \;, \label{eq smooth scheme}
\\
\ud \tilde{E}^{k,\epsilon}_t &=  \mu(\tilde{Y}^{k,\epsilon}_t,p) \ud t + \epsilon \ud B_t \label{eq smooth scheme E}
\end{align} 
with terminal condition $\tilde{Y}^{k,\epsilon}_h = \bar{v}^{k,\epsilon}(0,p,\tilde{E}^{k,\epsilon}_h)$ and initial condition 
$\tilde{E}^{k,\epsilon}_0:=e$.
\\
\noindent Observe that the $P$-variable is frozen in the above definition and that $\tilde{Y}^{k,\epsilon}_t = \tilde{v}^{k,\epsilon}(t,p,\tilde{E}^{k,\epsilon}_t)$, for $0 \le t \le h$.
\end{enumerate}

\noindent 
Before studying the truncation error, we give a strong local error control between the smooth approximations  $v^{k,\epsilon}$ and $\tilde{v}^{k,\epsilon}$.  Note that this local error control in $\sqrt{h}$ does not allow obtaining a converging global error control. We will however use it to obtain a better local control error in $L1$-norm, see the proof of Proposition 
\ref{pr truncation error}.
\begin{Lemma} \label{le minimal strong error on v-tilde v}
Under our \emph{standing assumptions} on $(\mu,b,\sigma)$, the following holds, for $p \in \R^d$, $h>0$,
\begin{align*}
\sup_{t\in [0,h], e \in \R}|v^{k,\epsilon}(t,p,e)-\tilde{v}^{k,\epsilon}(t,p,e)| \le C_{L_\psi}(1+|p|)\sqrt{h}.
\end{align*}
Importantly, $C_{L_\psi}$ does not depend on $k$ nor $\epsilon$, however it depends on the Lipschiz constant of $\psi$ in the $P$-variable.
\end{Lemma}

\proof 
For $t \le h$, let $V^{k,\epsilon}_t = \tilde{v}^{k,\epsilon}(t,p,\tilde{E}^{k,\epsilon}_t) - v^{k,\epsilon}(t,P^{\epsilon}_t,\tilde{E}^{k,\epsilon}_t) $ with $(P^\epsilon_0,\tilde{E}^{k,\epsilon}_0)=(p,e)$. 
Applying Ito's formula, we compute, since $\tilde{Y}^{k,\epsilon}$ is a martingale, recall \eqref{eq smooth scheme},
\begin{align*}
V^{k,\epsilon}_t = V^{k,\epsilon}_0 - \int_0^t\left(\partial_tv^{k,\epsilon}(s,P^{\epsilon}_s,\tilde{E}^{k,\epsilon}_s) + \mu(\tilde{Y}^{k,\epsilon}_s,p)\partial_e v^{k,\epsilon}(s,P^\epsilon_s,\tilde{E}^{k,\epsilon}_s) + \cL^\epsilon v^{k,\epsilon}(s,P^{\epsilon}_s,\tilde{E}^{k,\epsilon}_s) \right)\ud s + M^{k,\epsilon}_t,
\end{align*}
where $M^{k,\epsilon}$ is a square-integrable martingale.\\
From the PDE \eqref{value function pde} satisfied by $v^{k,\epsilon}$, we get
\begin{align}\label{eq dyn V}
V^{k,\epsilon}_t = V^{k,\epsilon}_0 - \int_0^t \left(
 \mu(\tilde{Y}^{k,\epsilon}_s,p)  - \mu(v^{k,\epsilon}(s,P^\epsilon_s,\tilde{E}^{k,\epsilon}_s),P^\epsilon_s)\right)\partial_e v^{k,\epsilon}(s,P^\epsilon_s,\tilde{E}^{k,\epsilon}_s)  \ud s + M^{k,\epsilon}_t\,.
\end{align}
We set, for $0 \le s \le h$, $ \delta P_s = P^\epsilon_s - p$ and
\begin{align}
c_s &:= \int_0^1\partial_y \mu(v^{k,\epsilon}(s,P^\epsilon_s,\tilde{E}^{k,\epsilon}_s) - \lambda V^{k,\epsilon}_s,p + \lambda \delta P_s) \ud \lambda \,,
\\
d_s &:= \int_0^1\partial_p \mu(v^{k,\epsilon}(s,P^\epsilon_s,\tilde{E}^{k,\epsilon}_s)- \lambda V^{k,\epsilon}_s,p + \lambda \delta P_s) \ud \lambda \;.
\end{align}
From  \eqref{eq conservation law}, we know that, for all $0 \le s \le h$,
\begin{align}\label{eq useful majo}
c_s \le -\ell_1< 0 \text{ and } |d_s|\le L\;.
\end{align}
Then, \eqref{eq dyn V} reads
\begin{align}
V^{k,\epsilon}_t = V^{k,\epsilon}_0 - \int_0^t 
 c_s V^{k,\epsilon}_s\partial_e v^{k,\epsilon}(s,P^\epsilon_s,\tilde{E}^{k,\epsilon}_s) \ud s 
 + \int_0^t d_s \delta P_s \partial_e v^{k,\epsilon}(s,P^\epsilon_s,\tilde{E}^{k,\epsilon}_s) \ud s
 + M^{k,\epsilon}_t
\end{align}
We set, for $0 \le t \le h$, $\cE_t = e^{\int_0^t c_s \partial_e v^{k,\epsilon}(s,P^\epsilon_s,\tilde{E}^{k,\epsilon}_s)  \ud s}$ and, we have
\begin{align}
0 \le \cE_t \le 1\;,\; \text{ for all } 0 \le t \le h\,,
\end{align}
since
$v^{k,\epsilon} \in \cK$, recall  \eqref{eq useful majo}.
We then compute
\begin{align}
\cE_t V^{k,\epsilon}_t = V^{k,\epsilon}_0 + \int_0^t d_s \delta P_s \partial_e v^{k,\epsilon}(s,P^\epsilon_s,\tilde{E}^{k,\epsilon}_s) \cE_s \ud s  + N^{k,\epsilon}_t
\end{align}
where $N^{k,\epsilon}$ is a square-integrable martingale. In particular, we get
\begin{align}\label{eq interm main}
|V^{k,\epsilon}_0| \le \esp{|V^{k,\epsilon}_h| - \int_0^h \frac{|d_s|}{|c_s|} |\delta P_s| c_s \partial_e v^{k,\epsilon}(s,P^\epsilon_s,\tilde{E}^{k,\epsilon}_s) \cE_s \ud s }
\end{align}
recall \eqref{eq useful majo}. Observe that, for all $0 \le s \le h$,
\begin{align}
\dot{\cE}_s := c_s \partial_e v^{k,\epsilon}(s,P^\epsilon_s,\tilde{E}^{k,\epsilon}_s) \cE_s
\end{align}
where the dot denotes classicaly the time derivative.
We thus deduce from \eqref{eq interm main}
\begin{align}
|V^{k,\epsilon}_0| &\le \esp{|V^{k,\epsilon}_h| + C \sup_{s\in[0,h]}|\delta P_s| (\cE_0-\cE_h) } \nonumber
\\
&\le \esp{|V^{k,\epsilon}_h|}  + C(1+|p|)\sqrt{h} \,. \label{eq almost there}
\end{align}
Now, we observe that
\begin{align}
\esp{|V^{k,\epsilon}_h|} &= \esp{|\bar{v}^{k,\epsilon}(0,p,\tilde{E}^{k,\epsilon}_h) - \psi(P^\epsilon_h,\tilde{E}^{k,\epsilon}_h)|}
\\
&=\esp{|\esp{\psi(P^\epsilon_h,\tilde{E}^{k,\epsilon}_h)} - \psi(P^\epsilon_h,\tilde{E}^{k,\epsilon}_h)|}
\\
&\le 2L_\psi \esp{|\delta P_h|}\;.
\end{align}
We thus get $\esp{|V^{k,\epsilon}_h|} \le C(1+|p|)\sqrt{h}$, which, combined with \eqref{eq almost there}, concludes the proof.
\eproof

\vspace{2mm}
\noindent
\textcolor{black}{We now turn to the main result for this part, which gives an upper bound to the error between $\cS_h(\psi)$ and $\Theta_h(\psi)$ that is effectively a control on the truncation error of the scheme.}
\begin{Proposition}[truncation error]\label{pr truncation error}
Under our standing assumptions on the coefficients $(\mu,b,\sigma)$, the following holds, for  $\psi \in \cK$:
\begin{align}\label{eq pr truncation}
\int |\cS_h(\psi)(p,e)-\Theta_h(\psi)(p,e)| \ud e \le C_{L_\psi}(1+|p|^2)h^{\frac32}\;,
\end{align}
for $p \in \R^d$, $h>0$.
\end{Proposition}

\proof 1. We first consider the regularised version of the decoupling fields, as introduced in \eqref{eq de bar v}-\eqref{eq smooth scheme}.
Let $V^{e,k,\epsilon}_t = \tilde{Y}^{k,\epsilon}_t - v^{\epsilon,k}(t,P^\epsilon_t,\tilde{E}^{k,\epsilon}_t)$, for $t \le h$, with $(P^\epsilon_t,\tilde{E}^{k,\epsilon}_t) =(p,e)$. We first observe that by definition \eqref{eq de bar v} and the fact that $P^\epsilon$ and $B$ (and thus subsequently $\tilde{E}^{e,k,\epsilon}$) are independent,
\begin{align}
\esp{V^{e,k,\epsilon}_h | B} &=  \bar{v}(0,p,\tilde{E}^{k,\epsilon}_h) - \esp{\psi^k(P^\epsilon_h,\tilde{E}^{k,\epsilon}_h) | B}
= 0 \label{eq term = 0} \;.
\end{align}
We also observe that, for $0 \le t \le h$,
\begin{align*}
|V^{e,k,\epsilon}_t| &\le  |\tilde{v}^{\epsilon,k}(t,p,\tilde{E}^{k,\epsilon}_t) - v^{\epsilon,k}(t,p,\tilde{E}^{k,\epsilon}_t)|
+
 |{v}^{\epsilon,k}(t,p,\tilde{E}^{k,\epsilon}_t) - v^{\epsilon,k}(t,P^\epsilon_t,\tilde{E}^{k,\epsilon}_t)|
 \\
 &\le C(1+|p|) \sqrt{h} + C|P^\epsilon_t-p|\;,
\end{align*}
where for the last inequality we used Lemma \ref{le minimal strong error on v-tilde v}  and the uniform Lipschitz continuity of $v^{\epsilon,k}$ (Note that $C$ depends upon the Lispchitz constant of $\psi$). This leads to
\begin{align}\label{eq key estim V}
\sup_{t \in [0,h]}\esp{\mathrm{esssup}_e|V^{e,k,\epsilon}_t|} \le C(1+|p|) \sqrt{h}\;.
\end{align}
2. Let us consider the tangent process $\partial_e \tilde{E}^{k,\epsilon}$ given by
\begin{align}
\partial_e \tilde{E}^{k,\epsilon}_t &= 1 + \int_0^t \partial_y \mu(\tilde{Y}^{k,\epsilon}_s,p)\partial_e\tilde{v}^{k,\epsilon}(s,p,\tilde{E}^{k,\epsilon}_s) \partial_e \tilde{E}^{k,\epsilon}_s \ud s 
\\
&= e^{\int_0^t\partial_y \mu(\tilde{Y}^{k,\epsilon}_s,p)\partial_e\tilde{v}^{k,\epsilon}(s,p,\tilde{E}^{k,\epsilon}_s) \ud s}.
\end{align}
And we observe that $0 \le \partial_e \tilde{E}^{k,\epsilon}_t \le 1$, for all $0 \le t \le h$. \\
In order to bound the error $\int |V^{e,k,\epsilon}_0| \ud e$, we will study the dynamics of
$t \mapsto \int |\esp{V^{e,k,\epsilon}_t\partial_e \tilde{E}^{k,\epsilon}_t} | \ud e$.
Using \eqref{eq dyn V}, we  compute
\begin{align}
V^{e,k,\epsilon}_t \partial_e \tilde{E}^{k,\epsilon}_t  &= V^{e,k,\epsilon}_0 + \int_0^t V^{e,k,\epsilon}_s\partial_y \mu(\tilde{Y}^{k,\epsilon}_s,p)\partial_e\tilde{v}^{k,\epsilon}(s,p,\tilde{E}^{k,\epsilon}_s) \partial_e \tilde{E}^{k,\epsilon}_s \ud s + N^{k,\epsilon}_t
\\
 & \!\!\!- \int_0^t \left( \mu(\tilde{Y}^{k,\epsilon}_s,p) -\mu(v^{k,\epsilon}(s,P^\epsilon_s,\tilde{E}^{k,\epsilon}_s),P^\epsilon_s)\right)
 \partial_e v^{k,\epsilon}(s,P^\epsilon_s,\tilde{E}^{k,\epsilon}_s)\partial_e \tilde{E}^{k,\epsilon}_s \ud s 
\end{align}
where $N^{k,\epsilon}$ is a square-integrable martingale. Taking expectation on both sides of the above equality,
we get
\begin{align}
|V^{e,k,\epsilon}_0| &\le
|\esp{ \int_0^h \left( \mu(\tilde{Y}^{k,\epsilon}_s,p) -\mu(v^{k,\epsilon}(s,P^\epsilon_s,\tilde{E}^{k,\epsilon}_s),P^\epsilon_s)\right)
 \partial_e v^{k,\epsilon}(s,P^\epsilon_s,\tilde{E}^{k,\epsilon}_s)\partial_e \tilde{E}^{k,\epsilon}_s \ud s }| 
 \\
 &+ |\esp{\int_0^h V^{e,k,\epsilon}_s\partial_y \mu(\tilde{Y}^{k,\epsilon}_s,p)\partial_e\tilde{v}^{k,\epsilon}(s,p,\tilde{E}^{k,\epsilon}_s) \partial_e \tilde{E}^{k,\epsilon}_s \ud s}| 
\end{align}
recall \eqref{eq term = 0}. Since $\partial_e\tilde{v}^{k,\epsilon}()$, $\partial_e \tilde{E}^{k,\epsilon}$ and $-\partial_y \mu()$ are non-negative,  we deduce
\begin{align}
|V^{e,k,\epsilon}_0| &\le
\esp{ \int_0^h  |\mu(\tilde{Y}^{k,\epsilon}_s,p) -\mu(v^{k,\epsilon}(s,P^\epsilon_s,\tilde{E}^{k,\epsilon}_s),P^\epsilon_s)|
 \partial_e v^{k,\epsilon}(s,P^\epsilon_s,\tilde{E}^{k,\epsilon}_s)\partial_e \tilde{E}^{k,\epsilon}_s \ud s } 
 \\
 &- \esp{\int_0^h |V^{e,k,\epsilon}_s|\partial_y \mu(\tilde{Y}^{k,\epsilon}_s,p)\partial_e\tilde{v}^{k,\epsilon}(s,p,\tilde{E}^{k,\epsilon}_s) \partial_e \tilde{E}^{k,\epsilon}_s \ud s}\,.
\end{align}
Integrating the previous inequality, we get
\begin{align}
\int |V^{e,k,\epsilon}_0| \ud e &\le
 \esp{ \int\int_0^h \left| \mu(\tilde{Y}^{k,\epsilon}_s,p) -\mu(v^{k,\epsilon}(s,P^\epsilon_s,\tilde{E}^{k,\epsilon}_s),P^\epsilon_s)\right|
 \partial_e [v^{k,\epsilon}(s,P^\epsilon_s,\tilde{E}^{k,\epsilon}_s)] \ud s \ud e }  =: A_1
 \nonumber
 \\
 &- 
\esp{ \int \int_0^h |V^{e,k,\epsilon}_s|\partial_e [ \mu(\tilde{v}^{k,\epsilon}(s,p,\tilde{E}^{k,\epsilon}_s),p)] \ud s\ud e} =: A_2
\,.
\end{align}
We know study the term $A_2$ above:
Since, $\tilde{v}^{\epsilon,k}$ is bounded and $\mu$ is Lipschitz continuous, we have
\begin{align}
|\int \partial_e [ \mu(\tilde{v}^{k,\epsilon}(s,p,\tilde{E}^{k,\epsilon}_s),p)] \ud e| \le C(1+|p|)
\end{align}
and then
\begin{align*}
A_2 &\le C(1+|p|)\esp{\int_0^h \mathrm{esssup}_e|V^{e,k,\epsilon}_s| \ud s}
\\
&\le C(1+|p|^2)h^{\frac32}\,,
\end{align*}
recalling \eqref{eq key estim V}.
\\
We now compute an upper bound  for $A_1$. Since $\mu$ is Lipschitz-continuous, we have
\begin{align*}
A_1  &\le C
 \esp{ \int\int_0^h \left( |V^{e,k,\epsilon}_s| + |P^\epsilon_s-p| \right)
 \partial_e [v^{k,\epsilon}(s,P^\epsilon_s,\tilde{E}^{k,\epsilon}_s)] \ud s \ud e } 
 \\
 &\le \esp{\int_0^h (\mathrm{esssup}_e|V^{e,k,\epsilon}_s| + |P^\epsilon_s-p|) \int \partial_e [v^{k,\epsilon}(s,P^\epsilon_s,\tilde{E}^{k,\epsilon}_s)] \ud e}
\end{align*} 
Since $v^{k,\epsilon}$ is (uniformly) bounded and using \eqref{eq key estim V}, we obtain
\begin{align}
A_1  &\le C(1+|p|^2)h^{\frac32}\;.
\end{align}
Combining the estimate for $A_1$ and $A_2$, we conclude:
\begin{align}
\int |\tilde{v}^{k,\epsilon}(0,p,e)-{v}^{k,\epsilon}(0,p,e)| \ud e \le C(1+|p|^2)h^{\frac32}\;.
\end{align}
Then passing to the limits in $k, \epsilon$ and using the dominated convergence theorem conclude the proof.
\eproof

\subsubsection{Scheme stability}

We now study the scheme's stability by a introducing a perturbed version of the scheme given in Definition \ref{de th splitting theo}.

\begin{Definition}[Perturbed  scheme] \label{de perturbed th splitting theo}
For $n \le N$, let $\eta_n: \R^d\times\R \rightarrow \R$ be measurable functions satisfying
\begin{align} \label{eq ass eta}
\int |\eta_n(p,e)|\ud e \le \mathfrak{c}(1+|p|^\kappa) \;\text{ for some } \kappa \ge 1, \mathfrak{c}>0,
\end{align}
where $\kappa$, $\mathfrak{c}$ do not depend on $n$.
We denote by $(u'_n)_{0 \le n \le N}$ the solution of the following backward induction:
\\
- for $n=N$, set $u'_N := \phi + \eta_N$,
\\
- for $n < N$, $u'_n = \cS_{t_{n+1}-t_n}(u'_{n+1}) + \eta_n$.
\end{Definition}

%

\begin{Proposition}[$L^1$-stability] \label{pr scheme stability}
Under our \emph{standing assumptions}, the following holds true for $(\eta_n)$,satisfying \eqref{eq ass eta}, perturbation of the scheme given in Definition \ref{de th splitting theo} :
\begin{align}
\max_{0 \le n \le N} \esp{\int |u^\pi_n - u'_n|(P^{0,p}_{t_n},e) \ud e} \le \sum_{n=0}^N \esp{\int |\eta_n|(P^{0,p}_{t_n},e)\ud e}\;.
\end{align}
\end{Proposition}
\proof 
1.a In the proof, we denote $\delta u^\pi_n = u_n - u'_n$, for all $n \le N$.\\
We observe that 
\begin{align} \label{eq bound on the scheme}
\int|\delta u_N|(p,e)\ud e &\le \int |\eta_N|(p,e)\ud e 
\le \mathfrak{c}(1 + |p|^\kappa)\,,
\end{align}
where we used \eqref{eq ass eta} for the last inequality. We have then
\begin{align}
\esp{\int|\delta u_N|(P^{0,p}_{T},e)\ud e} \le\esp{\int|\eta_N|(P^{0,p}_{T},e)\ud e} <\infty\;.
\end{align}
We used the fact that for any $q>0$,
\begin{align}\label{eq classic bound P}
\esp{\sup_{t \in [0,T]}|P^{0,p}_{t}|^q} \le C_q (1 + |p|^q).
\end{align}
\\
1.b Assume (induction hypothesis)
\begin{align}\label{eq induction hypothesis}
|\delta u_{n+1}|(p,e) \ud e \le K_{n+1}(1 + |p|^\kappa)\;.
\end{align}
for some positive $K_{n+1} < +\infty$ (Note that $K_N = \mathfrak{c}$ from \eqref{eq bound on the scheme}).
 Denoting $\bar{u}_{n+1} = \cD_{(t_{n+1}-t_n)}(u^\pi_{n+1})$ and $\bar{u}'_{n+1} = \cD_{(t_{n+1}-t_n)}(u'_{n+1})$, we have
\begin{align}
|u^\pi_n - u'_n | \le |\cT_{(t_{n+1}-t_n)}(\bar{u}_{n+1}) - \cT_{(t_{n+1}-t_n)}(\bar{u}'_{n+1})| + |\eta_n|
\end{align}
From Lemma 3.6 in \cite{jean2020modelling} applied to $\cT$, we obtain
\begin{align}
\int  |\cT_{(t_{n+1}-t_n)}(\bar{u}_n) - \cT_{(t_{n+1}-t_n)}(\bar{u}'_n)|(p,e) \ud e \le 
\int  |\bar{u}_{n+1} - \bar{u}'_{n+1}|(p,e) \ud e
\end{align}
Moreover,
\begin{align}
\bar{u}_{n+1} - \bar{u}'_{n+1} = \esp{u^\pi_{n+1}(P^{0,p}_{t_{n+1}-t_n},e)-u'_{n+1}(P^{0,p}_{t_{n+1}-t_n},e)}
\end{align}
which leads to 
\begin{align} \label{eq main step induction}
\int |\delta u_n|(p,e) \ud e \le \int\esp{|\delta u_{n+1}|(P^{0,p}_{t_{n+1}-t_n},e) \ud e} + \int |\eta_n|(p,e) \ud e\;.
\end{align}
From the induction hypothesis, we know that
\begin{align*}
\esp{|\delta u_{n+1}|(P^{0,p}_{t_{n+1}-t_n},e) \ud e} \le K_{n+1}(1 + \esp{|P^{0,p}_{t_{n+1}-t_n}|^\kappa})
\end{align*}
leading to $\esp{|\delta u_{n+1}|(P^{0,p}_{t_{n+1}-t_n},e) \ud e} \le K_{n+1} C_\kappa(1 + |p|^\kappa)$, 
where we used \eqref{eq classic bound P}. Using \eqref{eq ass eta}, 
we then obtain
\begin{align*}
\int |\delta u_n|(p,e) \ud e \le K_n(1 + |p|^\kappa) \text{ with } K_n = K_{n+1} C_\kappa + \mathfrak{c}\,,
\end{align*}
proving the induction hypothesis \eqref{eq induction hypothesis} for the next step. Moreover, this shows that
$\esp{\int |\delta u_{n+1}|(P^{0,p}_{t_{n+1}},e) \ud e} < +\infty$ and $\esp{\int |\delta u_{n}|(P^{0,p}_{t_{n}},e) \ud e} < +\infty$. Thus, we deduce from \eqref{eq main step induction},
\begin{align}\label{eq outcome induction}
\esp{|\delta u_{n}|(P^{0,p}_{t_{n}},e) \ud e} \le \esp{|\delta u_{n+1}|(P^{0,p}_{t_{n+1}},e) \ud e} + \esp{\int |\eta_n|(P^{0,p}_{t_{n}},e) \ud e}
\end{align}
2. From step 1.a and 1.b above, we deduce that \eqref{eq outcome induction} holds for all $n<N$. Iterating the inequality on $n$ concludes the proof.
\eproof

\subsubsection{Proof of Theorem \ref{th conv res theo scheme}}
We classicaly writes the true solution given by the decoupling field as a perturbed splitting scheme, for a perturbation $(\zeta_n)_{0 \le n \le N}$ given as follows: $\zeta_N = 0$ 
and $\zeta_n(\cdot) =  \Theta_{t_{n+1}-t_n}(\cV(t_{n+1},\cdot)) - \cS_{t_{n+1}-t_n}(\cV(t_{n+1},\cdot))$. We observe that, indeed, for all $n \le N$,
\begin{align}
\cV(t_n,\cdot) = \cS_{t_{n+1}-t_n}(\cV(t_{n+1},\cdot)) + \zeta_n(\cdot)\;.
\end{align}
From  Proposition \ref{pr truncation error}, we know that $\zeta_n$ satisfies  \eqref{eq ass eta} with $\kappa = 2$, recall \eqref{eq pr truncation}, and then
\begin{align}\label{eq trunc applied}
\esp{\int |\zeta_n|(P_{t_n},e) \ud e} \le C (1+|p|^2) (t_{n+1}-t_n)^{\frac32} .
\end{align}
Using Proposition \eqref{pr scheme stability}, we obtain that for $n=0$ in particular,
\begin{align*}
\int|\cV(0,p,e)-u^\pi_0(p,e)|\ud e \le C T (1+|p|^2)\sqrt{|\pi|}\,.
\end{align*}
\eproof

\section{Numerical schemes}
\label{se num scheme}

The possible difference in the dimension between the $E$-variable and the $P$-variable leads us to treat these variables very differently in the numerical procedure. The convergence result obtained in the previous section indicates that it is indeed reasonable to use a splitting scheme. We then work toward a fully implementable scheme building on this approach.

\subsection{A regression method for the splitting scheme}

We present here a -- still theoretical --  discrete-time scheme which combines a finite difference approximation of the transport operator and a probabilistic approximation of the diffusion operator. In the next section, we discuss various possible implementations. 
\vspace{1mm}

We first suppose that the approximation of the transport operator is given as follows.
Let $J$  be a positive integer and $\mathfrak{E}=(e_j)_{1 \le j \le J}$ a discrete grid of $\R$. We denote by $\cT^\mathfrak{E}_h$ an approximation of the operator $\cT_h$ on $\mathfrak{E}$. Namely,
\begin{align}
\R^d\times \R^J \ni (p,\theta) \mapsto \cT^\mathfrak{E}_h (p,\theta) \in \R^J.
\end{align}
This means that for each $p \in \R^d$, $\cT^\mathfrak{E}_h (p,\cdot)$ is an approximation on the grid $\mathfrak{E}$ of the corresponding equation \eqref{eq de cons law associated with transport} on $[0,h]$. We assume moreover that it satisfies, for some $q \ge 1$ and $q' \ge 1$,
\begin{align}\label{eq ass for wellposedness}
|\cT^\mathfrak{E}_h (p,\theta)| \le C(1 + |p|^q + |\theta|^{q'}) \;. 
\end{align}
The terminal condition $\psi:\R^d\times\R \rightarrow \R$ is simply approximated on $\mathfrak{E}$ by 
$\theta^j=\psi(p,e_j)$, for all $1 \le j \le J$ and $p \in \R^d$.

\vspace{2mm}
\noindent Given this approximate transport operator, we now introduce a probabilistic approximation of $\cV(0,p,\cdot)$ on $\mathfrak{E}$. To this end, let us consider the Euler scheme associated to $P$ on $\pi$, namely, for $n \ge 0$,
\begin{align}\label{eq euler scheme}
\widehat{P}^{\pi}_{t_{n+1}} = \widehat{P}^{\pi}_{t_{n}} + b(\widehat{P}^{\pi}_{t_{n}})(t_{n+1}-t_n)+ \sigma(\widehat{P}^{\pi}_{t_{n}}) \Delta \widehat{W}_n \; \text{ and } \widehat{P}^\pi_0 = p\,.
\end{align}
Here, $(\Delta \widehat{W}_n)_{0 \le n \le N-1}$ are independent random variables that stands for an approximation of the law of $(W_{t_{n+1}} - W_{t_n})_{0 \le n \le N-1}$ and we assume that their moments verify $\esp{|\Delta \widehat{W}_n|^\rho}\le C_\rho|t_{n+1}-t_n|^\frac{\rho}2$, $\rho \ge 1$. It is well known from the Lipschitz continuity assumption on $b$ and $\sigma$ that, for any $\rho \ge 1$,
\begin{align} \label{eq bound on P pi}
\esp{\sup_{t \in \pi}|\widehat{P}^\pi_t|^\rho }\le C_\rho(1+|p|^\rho)\;. 
\end{align}

\vspace{2mm}
\noindent We now define a discrete time process $(\Gamma_n)_{0 \le n \le N}$ valued in $\R^J$ as follows.
\begin{Definition}\label{de Gamma}
$(\Gamma_n)_{0 \le n \le N}$ is solution to the following backward scheme:
\begin{enumerate}
\item For $n=N$, $\Gamma^j_N=\phi(\widehat{P}^\pi_{t_N},e_j)$ for $1 \le j \le J$.
\item For $n < N$, compute
\begin{align}\label{splitting scheme}
{\Gamma}^j_n &= \cT^\mathfrak{E}_h(\widehat{P}^\pi_{t_n},\esp{\Gamma_{n+1} | \widehat{P}^\pi_{t_n}} )\;.
\end{align}
\end{enumerate}
\end{Definition}


\vspace{2mm}
\noindent For later use, we define the auxiliary process $(\bar{\Gamma}_n)$ by 
\begin{align}\label{eq step cond exp}
\bar{\Gamma}^j_n &= \esp{\Gamma^j_{n+1} | \widehat{P}^\pi_{t_n}} \;\text{ for all } 1 \le j\le J\,.
\end{align}

\noindent We also importantly observe that, due to the Markovian property of $\widehat{P}^\pi$ on $\pi$, $(\Gamma_n)$ satisfies  
\begin{align}\label{eq gamma func}
\Gamma_n := \gamma_n(\widehat{P}^\pi_{t_n})\;,\quad 0\le n \le N\,,
\end{align}
where the functions $ \gamma_n:\R^d \rightarrow \R^J$, are given by
\begin{Definition}\label{de functional scheme}
\begin{enumerate}
\item For $n=N$, $\gamma_N^j(p)=\phi(p,e_j)$, $1 \le j \le J$, $p \in \R^d$.
\item Then, compute for $n < N$, $p \in \R^d$,
\begin{align}
\bar{\gamma}^j_n(p) &= \esp{{\gamma}^j_{n+1}\left(p + b(p)h + \sigma(p) \Delta \widehat{W}_n \right)} \;\text{ for all } 1 \le j\le J,
\\
{\gamma}^j_n(p) &= \cT^\mathfrak{E}_h(p, \bar{\gamma}_n(p))\;.
\end{align}
\end{enumerate}
\end{Definition}

\noindent With the above definitions, we have that $\Gamma_0 = \gamma_0(P_0)$ which stands for an approximation of $\cV(0,P_0,\cdot)$ on the grid $\mathfrak{E}$.\\

\noindent To obtain the wellposedness of the previous definitions, we check that the conditional expectations at each step of the scheme are well defined. This follows from a direct backward induction using \eqref{eq ass for wellposedness} and \eqref{eq bound on P pi}.

\noindent 
Depending on how large $d$, the dimension of the P-variable, is, we may choose various probabilistic schemes to compute \eqref{eq step cond exp}. This has been thoroughly studied in the context of BSDEs approximation and various methods have been suggested: linear regression \cite{gobet2005regression,gobet2016stratified,gobet2016approximation}, quantization methods \cite{bally2003,BALLY20031,PAGES2018847}, cubature methods \cite{crisan2012solving,crisan2014second,chassagneux2020cubature} or Malliavin calculus approach \cite{BOUCHARD2004175,CRISAN}. In the next section, we present a non-linear regression method used e.g. in \cite{comehuredeepbsde}.


\subsection{Implementation using non linear regression}

We now turn to an implementation that can work in a high dimensional setting for $P$. To perform the regression step  in Definition \ref{de Gamma}, we will use Neural Networks representation of the value function. This will be coupled with conservative finite difference approximation of the transport operator that we first recall.

\subsubsection{Conservative Finite Difference approximation of transport equation}
We shall now consider conservative methods for the transport operator associated to the backward equation \eqref{eq de cons law associated with transport}. 

\vspace{2mm}
Recall that, for a given positive integer $J$, $\mathfrak{E}=(e_j)_{1 \le j \le J}$ is a uniform grid of $\R$ where we set $\delta  := e_{j+1}-e_j$. We also introduce 
$\Re = \set{r_0 = 0 <\dots<r_k<\dots<r_K=h}$ a uniform grid for a given positive integer $K$, and we set $\mathfrak{d} := h/K$.

\vspace{2mm}
\noindent We first consider the Lax-Friedrichs approximation to the backward transport equation \eqref{eq de cons law associated with transport} and define $\cT^{\mathrm{LF}}_{\mathfrak{E},\mathfrak{\Re},h}:\R^d\times\R^J \mapsto \R^J$ the approximation of the associated operator $\cT_h$. It is defined as follows, see e.g. \cite[Chapter 12]{leveque1992numerical}.

\begin{Definition}[Lax-Friedrichs] \label{de lf scheme}
For a given $p \in \R^d$ and $\theta \in \R^J$, let $(V^k_{j})_{1 \le k \le K, 1\le j \le J}$ denotes the approximation at the point $(r_k,e_j)$ The steps to compute $V$ are: 
\\
- at time $r_K=h$: set $V^K_{j} = \theta_j$, $1 \le j \le J$,
\\
- for $0\le k<K$: set $V^{k}_{1} = V^{k+1}_{1}$, $V^{k}_{J} = V^{k+1}_{J}$ and compute, for $1 < j < J$ : 
\begin{align}
V^{k}_{j} = \frac12(V^{k+1}_{j+1}+V^{k+1}_{j-1} )  
	+ \frac{\delta}{2 \mathfrak{d}}\left(\mathfrak{M}(p,V^{k+1}_{j+1}) - \mathfrak{M}(p,V^{k+1}_{j-1}) \right).
\end{align}
Then, set $\cT^{\mathrm{LF}}_{\mathfrak{E},\mathfrak{\Re},h}(p,\theta) := V^0$.
\end{Definition}

When the function $\mu$ has constant sign, a more satisfactory method to use is the upwind method, as it is less diffusive. Since in the application to carbon markets given in Example \ref{ex gbm P} and \ref{ex bm pos emit} below, this will be the case, we consider the upwind method for $\mu \ge 0$.
We thus now define $\cT^{\mathrm{U}}_{\mathfrak{E},\mathfrak{\Re},h}:\R^d\times\R^J \mapsto \R^J$ the approximation of the associated operator $\cT_h$ as follows, see again e.g. \cite{leveque1992numerical}.
\begin{Definition}[Upwind for $\mu \ge 0$]\label{de upwind scheme}
For a given $p \in \R^d$ and $\theta \in \R^J$ let $(V^k_{j})_{1\le j \le J, 1 \le k \le K}$ denotes the approximation at the point $(r_k,e_j)$ The steps to compute $V$ are: 
\\
- at time $r_K=h$: set $V^K_{j} = \theta_j$, $1 \le j \le J$,
\\
- for $0\le k<K$: set $V^{k}_{J} = V^{k+1}_{J}$ and compute, for $1 \le j < J$ : 
\begin{align}
V^{k}_{j} = V^{k+1}_{j} 
	+ \frac{\delta}{ \mathfrak{d}}\left(\mathfrak{M}(p,V^{k+1}_{j+1}) - \mathfrak{M}(p,V^{k+1}_{j}) \right).
\end{align}
Then, set $\cT^{\mathrm{U}}_{\mathfrak{E},\mathfrak{\Re},h}(p,\theta) := V^0$.
\end{Definition}


\subsubsection{Non-linear regression and implemented scheme}
We first mention that for this part the Euler scheme \eqref{eq euler scheme} is computed using real Brownian increment, namely $\widehat{\Delta W}_n = (W_{t_{n+1}}-W_{t_n})$, $0 \le n \le N-1$. We have seen in the last section two possible implementations of the transport operator $\cT_h$ on the spatial grid $\mathfrak{E}$, that we shall denote for this part simply by $\cT^\mathfrak{E}_h$.
The last point to precise is the computation of the conditional expectation part of the scheme in Definition \ref{de Gamma}, 
where at each time step the quantities $\gamma_n(\widehat{P}^\pi_{t_n})=\esp{\Gamma_{n+1} | \widehat{P}^\pi_{t_n}}$ has to be estimated, recall Definition \ref{de functional scheme}.
In order to do so, we will use deep learning as it was demonstrated to be very efficient for high dimensional system, already in the setting of FBSDEs, see e.g. \cite{han2018solving, comehuredeepbsde}. The functions $(\gamma_n)$ will be optimally approximated by a feedforward neural network. We denote by $\cN\cN_{d_0,d_1,L,m}$ the set of neural nets, which are functions $\Phi(\cdot;\Theta):\R^{d_0} \mapsto \R^{d_1}$, parametrised by $\Theta$ and  with the following characteristics: \textcolor{black}{the input dimension is $d_0$, the output dimension is $d_1$, $L+1$ is the number of layers, $m=(m_l)_{0\le m_l\le L}$ where $m_l$ is the number of neurons on each layer, $l=0,...,L$: by default, $m_0 = d$ and $m_L = d_1$. The neural network has thus $L-1$ hidden layers. We refer, to e.g. \cite[Section 2]{comehuredeepbsde} for a detailed description of feedforward neural network. The number of total parameters is $N_{L,m} =\sum_{l=0}^{L-1}m_l(1+m_{l+1})$, and thus $\Theta \in \R^{N_{L,m}}$.} 

\vspace{2mm}
Given $\cT^\mathfrak{E}_h = \cT^{\mathrm{U}}_{\mathfrak{E},\mathfrak{\Re},h}$ or $\cT^\mathfrak{E}_h = \cT^{\mathrm{LF}}_{\mathfrak{E},\mathfrak{\Re},h}$, the scheme to compute $(\hat{\gamma}_n,\hat{\bar{\gamma}}_n$  approximation of $(\gamma_n,\bar{\gamma}_n)$ in Definition \ref{de functional scheme} is given as follows.

\begin{Definition}\label{de deep algo} 
\begin{enumerate}
	\item At $t_N= T$, $\widehat{\bar{\gamma}}^j_n (p) = \widehat{\gamma}_N^j(p)=\phi(p,e_j)$, $1 \le j \le J$, $p \in \R^d$.
	\textcolor{black}{
	\item For $n=N-1,...,1$: given a simulation of $P^\pi_{t_n}$, optimize 
	\begin{align}
		\widehat{\mathcal{L}}_n(\Theta) = \mathbb{E}\Big|
		\cT^\mathfrak{E}_h(\widehat{P}^\pi_{t_{n+1}}, \widehat{\bar{\gamma}}_{n+1}(\widehat{P}^\pi_{t_{n+1}})) - 
		 \left( \cY_{n}(\widehat{P}^\pi_{t_{n}},\Theta) + \cZ_n(\widehat{P}^\pi_{t_{n}},\Theta) (W_{t_{n+1}}-W_{t_n})\right)
		\Big|^2 \label{eq loss}
  \end{align}
		where 
		$(\cY_{n}(\cdot,\Theta),\cZ_{n}(\cdot,\Theta)) \in \cN\cN_{d,(d+1)\times J,L,m}$, 
  so that
  \begin{align*}
		\Theta_n^{\star}\in \arg \min_{\Theta\in\mathbb{R}^{N_m}}\widehat{\mathcal{L}}_n(\Theta) 
		\quad\text { and  then } \quad \widehat{\bar{\gamma}}_n(\cdot) := \cY_{n}(\cdot,\Theta_n^{\star})\,.
	\end{align*}
	}	
	\item At the initial time $t_0=0$, compute $\widehat{\gamma}_0(\widehat{P}^\pi_0) = \esp{\cT^\mathfrak{E}_h(\widehat{P}^\pi_{t_{1}}, \widehat{\bar{\gamma}}_{1}(\widehat{P}^\pi_{t_{1}}))}$.
	 \end{enumerate}
\end{Definition}
The function $\hat{\gamma}_0(\cdot)$ stands for the numerical approximation of $\cV(0,P_0,\cdot)$.

%
%

{\color{black}
\begin{Remark}
\begin{enumerate}[(i)]
\item The loss minimisation in \eqref{eq loss} is done using a Stochastic Gradient Descent algorithm: we use Adam Optimizer \cite{kingma2014adam} provided in the Keras API \cite{chollet2015keras}. The good approximation of the function $\bar{\gamma}_n$ is guaranteed by universal approximation theorem for neural networks \cite{hornik1989multilayer} and is quite efficient in practice, as demonstrated by our numerical examples below.
\item Definition \ref{de deep algo} should be compared with the scheme DBDP1 in \cite[Section 3]{comehuredeepbsde}. In this paper, the authors compute a non-linear conditional expectation (related to BSDEs) at each step: Here, we only compute a conditional expectation and use $\cZ$ as a control variate. In particular, differently to DBDP1, we have to apply $\cT^\mathfrak{E}_h$ in $\widehat{\cL}_n(\cdot)$ at each step. Note also that one can not apply directly DBDP1 to the singular FBSDE \eqref{eq singular fbsde} under study as it is a fully-coupled FBSDE.
\item A key point is to ensure the stability of the finite difference scheme for the transport equation, namely that the CFL condition is satisfied. For example, for the Lax-Friedrichs scheme given in Definition \ref{de lf scheme}, one has to enforce, for each time $t_n$:
	 \begin{align*}
	 	\sup_{1\leq k\leq K, 1\leq j\leq J}\Big|\mu(V_{k}^j,\widehat{P}^\pi_{t_{n}})\frac{\delta}{\mathfrak{d}}\Big|<1,
	 \end{align*}
see e.g. \cite[Chapter 13]{leveque1992numerical}.
In practice, we choose $B$, such that
\begin{align*}
\sup_{y \in [0,1], p \in [-B,B]^d}|\mu(y,p)\frac{\delta}{\mathfrak{d}}\Big|<1.
\end{align*} 
The constant $B$ depends obviously on the parameters $\delta, \mathfrak{d}$ and should be large enough. Then, in the simulation, $\widehat{P}^\pi_{t_n}$ is projected on $[-B,B]^d$. We also ensure that $0\le V_{k}^j \le 1$ by truncating $\bar{\gamma}^j_n$ if necessary and relying on the scheme monotony.
\end{enumerate}
\end{Remark}
}

\subsection{Numerical experiments}
\label{se num exp}

In this section, we present the results of our numerical experiments that show that the splitting scheme is efficient in practice to approximate $\cV(\cdot)$.
The method presented in the previous section, will be tested on two complementary models to  Example \ref{ex linear emission}. The first one reads as follows.
\begin{Example}[BM with positive emission]\label{ex bm pos emission}
\label{ex bm pos emit}
\begin{align}\label{eq de mupos w any d}
\ud P^\ell_t = \sigma \ud W^\ell_t \text{ and } \ud E_t = \mu(\frac1{\sqrt{d}} \sum_{\ell=1}^d P^\ell_t,Y_t) \ud t
\end{align}
with $\mu(p,y) = 1 + \frac1{1+e^{-p}} - y$ and $\phi(p,e) = \1_{\set{e \ge 0}}$\;.
\end{Example}
The above model will have non negative $\mu$ which is more realistic if one has in mind application to carbon market. A critic could be however that it is driven by a Brownian Motion and that it will not suffer any discrete time error. We then introduce a multiplicative model as follows.
\begin{Example}[Multiplicative model]\label{ex mulitplicative model}
\label{ex gbm P}
\begin{align}
\ud P^\ell_t = \mu P^\ell_t \ud t + \sigma P^\ell_t \ud W^\ell_t, \,P^\ell_0=1, \text{ and } \ud E_t = \tilde \mu(P_t,Y_t) \ud t
\end{align}
with
$
\tilde{\mu}(p,y)= \left(\prod_{\ell=1}^d p^\ell \right)^{\frac1{\sqrt{d}}}e^{-\theta y}
$, for some $\theta > 0$ and $\phi(p,e) = \1_{\set{e \ge 0}}$\;.
\end{Example} 
We are not aware of any explicit solution for these models but they have the property that any $d+1>2$ dimensional model can be recast as a $2$-dimensional model (one dimension for the $P$-variable, one dimension for the $E$-variable). This will be used  for numerical validation of the non-linear regression scheme used for the multidimensional models by introducing an alternative scheme efficient in low dimension (see next Section). However, one should notice that there is no simple equivalent one-dimensional PDE available for Examples \ref{ex gbm P} or \ref{eq de mupos w any d} as it is the case for Example \ref{ex lin model}, recall \eqref{eq pde change variable}.

\subsubsection{An alternative scheme}
\label{subse alternative scheme}
To validate empirically the results obtained with the non-linear regression scheme, we could use a PDE method in low dimension. However, we chose to use here another method based on the splitting scheme that will combine a particle method with tree-like regression. This method will be efficiently implemented on the Examples \ref{ex bm pos emission}, \ref{ex mulitplicative model} and \ref{ex lin model} for two main reasons: We  work in moderate dimension and the process $P$ can be expressed as a function of the underlying Brownian motion, namely $P_t = \mathfrak{P}(t,W_t).$
\\
Here, contrarily to the previous section,  $(\Delta \widehat{W}_n := \widehat{W}_{t_{n+1}} - \widehat{W}_{t_n})_{0 \le n \le N-1}$ stands for discrete approximation of the  Brownian increments $(W_{t_{n+1}} - W_{t_n})_{0 \le n \le N-1}$. 
We also assume that the time grid $\pi$ is equidistant and thus $|\pi|=\frac{T}N=:\mathfrak{h}$.
One could then use,
for all $1\leq \ell \leq d$,  $\mathbb{P}(\widehat{W}_n^{\ell}= \sqrt{h}) = \mathbb{P}(\widehat{W}_n^{\ell}= -\sqrt{h}) = \frac{1}{2}$ in \eqref{eq euler scheme} but this requires $2^d$ points in total for the approximation. We use instead the cubature formula introduced in \cite[Section A.2]{gyurko2011efficient} which requires only $2d$ points. Denoting $(\mathfrak{e}^\ell)_{1\le \ell \le d}$ the canonical basis of $\R^d$, we set, for $1 \le i \le I=2d$,
$\P(\Delta \widehat{W}_n=\omega^i_{\mathfrak{h}}) = \frac1{2d}$ and $\omega^i_{\mathfrak{h}} = -\sqrt{d{\mathfrak{h}}}\mathfrak{e}^\ell$, if $i=2\ell$ or $\omega^i_{\mathfrak{h}} = \sqrt{d{\mathfrak{h}}}\mathfrak{e}^\ell$ if $i=2\ell - 1$. At a given point $(t_n,\widehat{P}^{\pi}_{t_{n}}=\mathfrak{P}(t_n,\widehat{W}_{t_n}),e)$, the approximation of $\cD_{\mathfrak{h}}(\psi)$, recall Definition \ref{de diffusion step}, reads then simply
\begin{align}\label{eq cubature estimation}
\esp{\psi(\widehat{P}^{\pi}_{t_{n+1}},e)|\widehat{P}^{\pi}_{t_{n}}} = 
\frac1{2d} \sum_{i=1}^I \psi(\mathfrak{P}(t_{n+1},\widehat{W}_{t_n}+\omega^i_{\mathfrak{h}}),e)\;.
\end{align}
For $0 \le n \le N$, we denote by ${\mathfrak{S}}_n$,  the discrete support of the random variable $\widehat{W}_{t_n}$.
We observe that $\mathfrak{S}_n \subset \mathfrak{S}_{n+1}$ and for $x \in \mathfrak{S}_n$, $x + \Delta \widehat{W}_n \in \mathfrak{S}_{n+1}$. Thus, when computing \eqref{eq cubature estimation}, there is no need for an interpolation step, if $\psi(\cdot,e)$ is known on $\mathfrak{S}_{n+1}$. We will obviously exploit this fact and compute recursively, backward in time, the approximation of $\cV$ on the discrete sets $({\mathfrak{S}}_n)_{0\le n\le N}$. The  full approximation of the diffusion operator $\cD_{(t_{n+1}-t_n)}$ that acts at time $t_n$, will be given after discussing the discrete version of the operator $\cT$, as it will be then more easily justified. 

\vspace{2mm} 
Let us thus now introduce a discrete version of the operator $\cT$, recall Definition \ref{de transport step}, that will compute an approximation to \eqref{eq de cons law associated with transport} written in \emph{forward form}: We shall use the celebrated Sticky Particle Dynamics (SPD) \cite{brenier1998sticky}  see also \cite[Section 1.1]{jourdain36optimal}. The SPD is particularly simple to implement in our case, since, due to the monotonicity assumption on $(\mu,\psi)$, there is no particle colliding!  For $M \ge 1$, let $D_M=\set{\mathrm{e}=(e_1,\dots,e_m,\dots,e_M) \in \R^M \,|\, e_1 \le \dots \le e_m \le \dots \le e_M}$. The discrete version of $\cT$ will act on empirical CDF or equivalently on empirical distribution $\frac1M\sum_{m=1}^M \delta_{e_m}$ ($\delta_e$ is the Dirac mass at $e$). Generally, $\psi(p,.)$, which is a CDF for each $p \in \R^d$, would need to be approximated in an optimal way on $D_M$. We observe here that the terminal condition $\phi$ to Examples \ref{ex bm pos emission}, \ref{ex mulitplicative model} and \ref{ex lin model}, is simply represented by $\mathrm{e}=(0,\dots,0)$. The iterative algorithm allows us to restrict our study to terminal condition $\psi$, such that $\psi(p,\cdot) = {H}\!*\!(\frac1M \sum_{m=1}^M \delta_{e_m})$ for some $\mathrm{e}\in D_M$, where $H$ is the Heaviside function and $*$ the convolution operator.
The approximation of $\cT$ is then given by
\begin{align}
\R^d \times D_M \ni (p,\mathrm{e}) \mapsto \cT_h^{M}(p,\mathrm{e})=(E^{p,m}_h)_{1 \le m \le M} \in D_M\;,
\end{align}
where $(E^{p,m}_h)_{1 \le m \le M}$ is a set (of positions) of particles computed as follows.
Given the initial position $\mathrm{e} \in D_M$ (representing $\psi$) and velocities $(\bar{F}_m(p))_{1 \le m \le M}$ set to 
$\bar{F}_m(p) = - \int_{(m-1)/M}^{m/M}\mu(p,y)\ud y$, we consider $M$ particles $(E^{p,m})_{1 \le m \le M}$, whose positions at time $t \in [0,h]$ are simply given by
\begin{align}\label{eq particle dynamics}
E^{p,m}_t = e_m + \bar{F}_m(p) t\;.
\end{align}
We observe that $(E^{p,m}_t)_{1 \le m \le M} \in D_M$, for all $t \in [0,h]$, as $-\mu$ is non-decreasing.

We are now ready to define  the approximation of the diffusion operator $\cD_{t_{n+1}-t_n}$, denoted $\cD^M_n$: it  will  take into account
that $\cT_h^M$ acts at the level of particles.  Introduce, to ease the presentation, $\cP_{n+1}=\set{p = \mathfrak{P}(t_{n+1},\mathrm{w}),\mathrm{w} \in \mathfrak{S}_{n+1}}$, which is simply the (discrete) support of $\widehat{P}^\pi_{t_{n+1}}$. Assume that 
$$\mathfrak{S}_{n+1} \ni \mathrm{w} \mapsto \Psi(\mathrm{w}) = \mathrm{e}^\mathrm{w} \in D_M$$ 
is given such that for $p \in \cP_{n+1}$, $p=\mathfrak{P}(t_{n+1},\mathrm{w})$, we have $\psi(p,\cdot) = {H}\!*\!(\frac1M \sum_{m=1}^M \delta_{e^{\mathrm{w}}_m})$. Then,  for  $\mathrm{w} \in \mathfrak{S}_n$, setting $\mathrm{w}^i_{n+1} = \mathrm{w}+\omega^i_{\mathfrak{h}}$, \eqref{eq cubature estimation} reads
\begin{align*}
\bar{v}_n(\mathrm{w},e) &:=\esp{\psi(\widehat{P}^{\pi}_{t_{n+1}},e)|\widehat{P}^{\pi}_{t_{n}} 
= \mathfrak{P}(t_{n},\mathrm{w})} 
\\ & = 
\frac1{2d} \sum_{i=1}^I H\!*\!(\frac1M \sum_{m=1}^M \delta_{e^{\mathrm{w}^i_{n+1}}_m})(e)\;,
\\
&=  H\!*\!(\frac1{2dM} \sum_{i=1}^I \sum_{m=1}^M \delta_{e^{\mathrm{w}^i_{n+1}}_m})(e)\;.
\end{align*}
This means that,  the function $e \mapsto \bar{v}(\mathrm{w},e)$ is an empirical CDF and is determined by the particles $\cE=\bigcup_{i=1}^I \set{\mathrm{e}^{\mathrm{w}^i_{n+1}}}$. There is no need to keep $2dM$ particles at step $n$, when the function $\psi$ at step $n+1$ is given by $M$ particles (for each $p \in \cP_{n+1}$). To reduce the number of particles, we first sort the cloud of particles $\cE$ to obtain $\tilde{\mathrm{e}}^{\mathrm{w}} \in D_{2dM}$, then we consider $\bar{\mathrm{e}}^{\mathrm{w}}:=(\tilde{e}^{\mathrm{w}}_{2dm})_{1 \le m \le M}$.
The approximation operator $\cD_n^M$ is finally defined by
\begin{align}\label{eq de particle regression definition}
(D_M)^{\mathfrak{S}_{n+1}} \ni \Psi \mapsto \cD_n^M(\Psi)(\mathrm{w})=\bar{\mathrm{e}}^{\mathrm{w}} \in 
(D_M)^{\mathfrak{S}_{n}}\;.
\end{align}

\vspace{2mm} 
\noindent The overall procedure is as follows
\begin{Definition}[Alternative scheme] \label{de art scheme}
\begin{enumerate}
\item  At $n=N$: Set $\mathrm{e}_N :=(0,\cdots,0)$ whose empirical CDF is $\phi$. 
Then define, $\gamma_N$ by
\begin{align*}
\mathfrak{S}_{N} \ni \mathrm{w} \mapsto \gamma_N(\mathrm{w}) 
= \mathrm{e}_N\;.
\end{align*}

\item For $n<N$: Given $\gamma_{n+1}:\mathfrak{S}_{n+1} \rightarrow D_M$, define
$\bar{\gamma}_n:\mathfrak{S}_{n} \rightarrow D_M$ by $\bar{\gamma}_n = \cD^M_n(\gamma_{n+1})$ and
then ${\gamma}_n$ by
\begin{align}
\mathfrak{S}_{n} \ni \mathrm{w} \mapsto \gamma_n (\mathrm{w}) 
	= \cT_{\mathfrak{h}}^M(\mathfrak{P}(t_n,\mathrm{w}),\bar{\gamma}_n (\mathrm{w})) \in D_M\;.
\end{align}
\end{enumerate}
The approximation of $\cV(0,0,\cdot)$ is then given by $H\!*\!\gamma_0$. 
 \end{Definition}
 
 \subsubsection{Numerical results}
 
 In this section, we report the findings of the numerical tests we performed on the models given in Examples \ref{ex bm pos emission}, \ref{ex mulitplicative model} and \ref{ex lin model}, using the 
 non-linear regression and splitting method of Definition \ref{de Gamma} and the alternative scheme, presented in Definition \ref{de art scheme}.
 
 \vspace{2mm}
 Concerning the non-linear regression, we use a common structure in all our experiments for the feedforward neural networks used in \eqref{eq loss} to represent $(\cY,\cZ)$, namely:  
\\
- The output layer is of dimension $(J+1)\times d$, where $J$ is size of the $E$-variable grid;
\\
- Two intermediate layers of dimension $\kappa_J\times d + 10$ ($\kappa_J$ is fixed to 20 below);
\\
- An input layer of dimension $d$.
\\
As already mentionned, the training is done using the Adam optimiser using 100 mini-batches with size 50 and  batch normalization. We check validation loss every 30 iterations with the  validation batch of size 500.
The learning rate is initially fixed at $\eta=0.001$. 
%

\begin{figure}[h]
\centering
\subfloat[$\sigma=0.01$]{\label{a}\includegraphics[width=.3\linewidth]{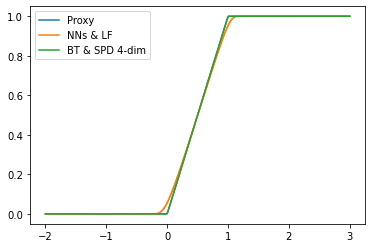}}  \hfill
\subfloat[$\sigma=0.3$]{\label{b}\includegraphics[width=.3\linewidth]{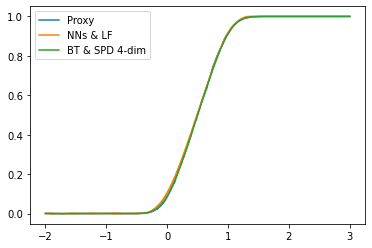}}\hfill 
\subfloat[$\sigma=1.0$]{\label{c}\includegraphics[width=.3\linewidth]{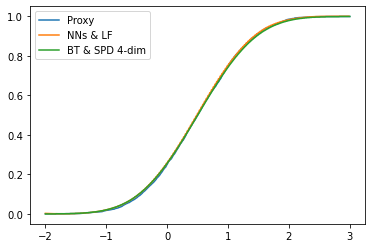}}
\caption{
Model of Example \ref{ex lin model}: Comparison of the two methods Neural Nets \& Lax-Friedrichs (NN\&LF) with $d=10$ and the alternative scheme (BT\&SPD) with $d=4$. The Proxy solution is given by the same particle method used in Figure \ref{fig main figure} on the one-dimensional PDE \eqref{eq pde change variable}. Lax-Friedrichs scheme implemented with discretization of space $J=1500,1000,500$, for $\sigma=0.01,0.3, 1$ respectively and number of time step $K=30$. The number of time step for the splitting is $N=64$. For $BT\&SPD$, the number of particles is $M=3500$ and the number of time steps $N= 20$.}
\label{fig neural nets d=10}
\end{figure}

\vspace{2mm}
We first observe that our schemes are able to reproduce the proxy for the true solution of Example \ref{ex lin model} as reported in Figure \ref{fig neural nets d=10}. This has to be compared with the results of Figure \ref{fig main figure} for the ``classical'' FBSDEs methods. Let us emphasize that the non-linear regression scheme (denoted $NN\&LF$) is tested in dimension $d=10$ and the alternative scheme (denoted $BT\&SPD$) in dimension $d=4$. Since the Lax-Friedrichs scheme presents a diffusive phenomenon, we increase the space discretization steps to overcome this effect when $\sigma$ decreases. In the numerical computations, we choose respectively $J=500, 1000, 1500$ for $\sigma=1.0, 0.3, 0.01$ to obtain satisfactory approximations.

\begin{figure}[h]
	\centering
	\subfloat[$\sigma=0.01$]{\label{a}\includegraphics[width=.3\linewidth]{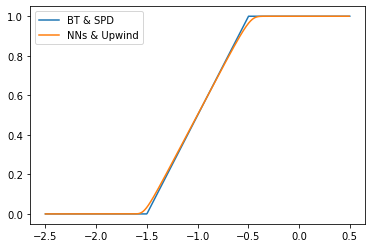}}\hfill
	\subfloat[$\sigma=0.3$]{\label{b}\includegraphics[width=.3\linewidth]{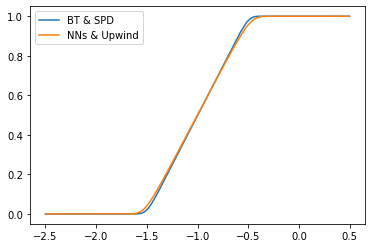}}\hfill 
	\subfloat[$\sigma=1.0$]{\label{c}\includegraphics[width=.3\linewidth]{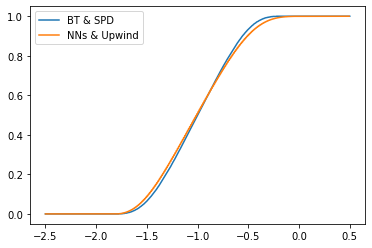}}
	\caption{Example \ref{ex bm pos emission} in dimension $d=10$: Comparison of two methods Neural nets \& Upwind scheme and  solution obtained using the alternative scheme on equivalent 4-dimensional model. The Upwind scheme used discretization of space $J=100,300,400$ respectively for $\sigma = 1,0.3,0.01$ and number of time step $K=20$. The number of time step for the splitting is $N=32$. For $BT\&SPD$, the number of particles is $M=3500$, and the number of time steps $N=20$.}	
	\label{fig bm pos emission}
\end{figure}
\begin{figure}[h]
	\centering
	\subfloat[$\sigma=0.01$]{\label{a}\includegraphics[width=.3\linewidth]{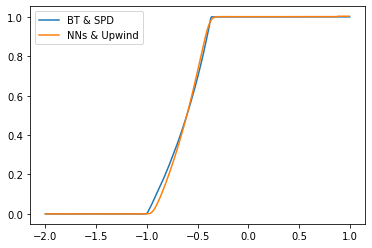}}\hfill
	\subfloat[$\sigma=0.3$]{\label{b}\includegraphics[width=.3\linewidth]{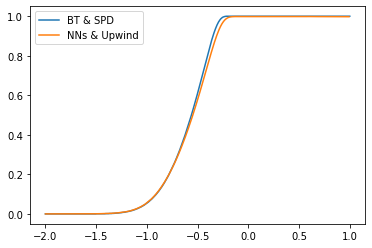}}\hfill 
	\subfloat[$\sigma=1.0$]{\label{c}\includegraphics[width=.3\linewidth]{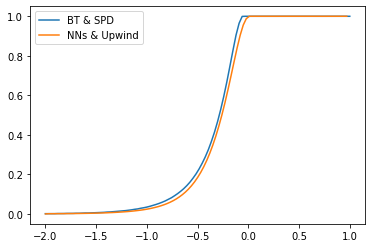}}
	\caption{Example \ref{ex mulitplicative model} in dimension $d=10$: Comparison of two methods Neural nets \& Upwind scheme and  solution obtained using the alternative scheme on equivalent 4-dimensional model (BT\&SPD). The Upwind scheme used discretization of space $J=100,400,500$ respectively for $\sigma = 1,0.3,0.01$ and number of time step $K=20$. The number of time step for the splitting is $N=32$. For $BT\&SPD$, the number of particles is $M=3500$, 
and the number of time steps $N= 20$.	
	}
	\label{fig multiplicative model}
\end{figure}

Next, we tested our scheme on the models of Example  \ref{ex bm pos emission} and  Example \ref{ex mulitplicative model}. The results are reported on the graphs in Figure \ref{fig bm pos emission} and \ref{fig multiplicative model} respectively.
Since the function $\mu$ is always positive in these two examples, we can use an Upwind scheme. Unlike  the Lax-Friedrichs scheme, the Upwind scheme is less diffusive, and we can lower the number of space discretization. In our example, we take $J=100,300,400$ respectively for $\sigma=1.0, 0.3, 0.01$. We are not aware of an exact solution  for this model, so we compare both the non-linear regression scheme (NN\&U) for $d=10$ and the alternative scheme (BT\&SPD) on an equivalent four dimensional model.

\begin{figure}[h]
	\centering
	\subfloat[$\sigma=0.01, J = 500$]{\label{c}\includegraphics[width=.3\linewidth]{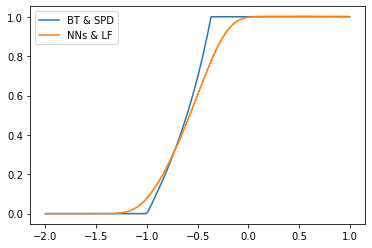}}
		\subfloat[$\sigma=0.01, J = 1500$]{\label{c}\includegraphics[width=.3\linewidth]{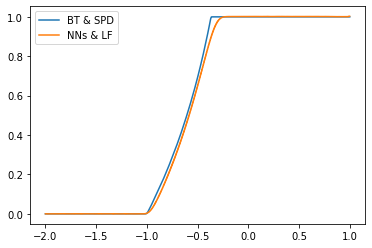}}
	\caption{Example \ref{ex mulitplicative model} in dimension $d=10$: Neural nets \& Lax-Friedrichs with $J=500$ and $1500$, compared with the Proxy (BT\&SPD in dimension one).}
	\label{fig illustration diffusive behavior}
\end{figure}

As we pointed out before, Lax-Friedrichs scheme is more diffusive than Upwind scheme: this is illustrated on Figure \ref{fig illustration diffusive behavior}, by considering the case where $\sigma=0.01$, and taking $J=400$ only for the $LF$ space discretisation. On this graph and the computations below, the `Proxy'  to the true solution is obtained by running an equivalent one-dimensional model using the alternative scheme (BT\&SPD) with parameters: number of particles $M=3500$, number of time step $N=20$. 
\textcolor{black}{Table \ref{the main table} presents the error obtained by comparing the non-linear regression scheme to this proxy, the computational times is also given\footnote{Intel Core i5-8265U, 16.0 GB RAM.} The $L1$-error is the error used in the theoretical part, but one can see that the $L\infty$ error behaves also very well. The computational times can still be reduced on our examples by diminishing the batch size but this would certainly not generalise to other models more challenging for the training of the neural networks.}

\begin{center}
	\begin{tabular}{|c| c| c| c| c| c| c|} 
		\hline
		Model & Sigma & Method & Parameters & $L_{1}$ & $L_{\infty}$ & Time \\ [0.5ex] 
		\hline\hline
		Ex \ref{ex lin model} &1.0 & NN $\&$ LF & $J=500$ & 0.0233
		& 0.0283 & 3813s\\
		\hline
		Ex \ref{ex bm pos emission} & & NN $\&$ Upwind & $J=500$ & 0.0116 & 0.0142 & 1687s\\
		\hline
		Ex \ref{ex mulitplicative model} & & NN $\&$ Upwind & $J=100$ & 0.0206
		& 0.0856 & 336s \\
		\hline
		
		Ex \ref{ex lin model} & 0.3 & NN $\&$ LF & $J=1000$ & 0.0183 & 0.0250 & 7660s \\
		\hline
		Ex \ref{ex bm pos emission} & & NN $\&$ Upwind & $J=500$ & 0.0147 & 0.0220 & 1693s \\
		\hline
		Ex \ref{ex mulitplicative model} & & NN $\&$ Upwind & $J=400$ & 0.0756 & 0.1253 & 1488s\\
		\hline
		
		Ex \ref{ex lin model} & 0.01 & NN $\&$ LF & $J=2000$ & 0.0055 & 0.0215 & 15139s\\
		\hline
		Ex \ref{ex bm pos emission} & & NN $\&$ Upwind & $J=500$ & 0.0141 & 0.0365 & 1712s \\
		\hline
		Ex \ref{ex mulitplicative model} & & NN $\&$ Upwind & $J=500$ & 0.0410 & 0.0843 & 1701s\\
		\hline
		
	\end{tabular}
	\captionof{table}{Numerics of model   \ref{ex lin model}, \ref{ex bm pos emission} and \ref{ex mulitplicative model} with different parameters}\label{the main table}
\end{center}

Finally, we want to empirically estimate the convergence rate of the error introduced by the splitting.
We consider the model \ref{ex mulitplicative model} where $\sigma=0.3$ and for which there is no discrete-time simulation error (as the forward process is a Brownian Motion). We consider a set of number of time steps $N:= \{4,8,16,32,64,128\}$, and compute the $L^{1}$ and $L^{\infty}$ error by NN $\&$ Upwind method (with $K=20$, $J=400$). The proxy solution is always given by alternative scheme in one dimensional equivalent model, to achieve a better precision. The empirical convergence rate with respect to the number of time step is close to one, which is slightly better than the upper bound obtained in Theorem \ref{th conv res theo scheme}.
\begin{figure}[h]
	\centering
	\includegraphics[width=0.3\textwidth]{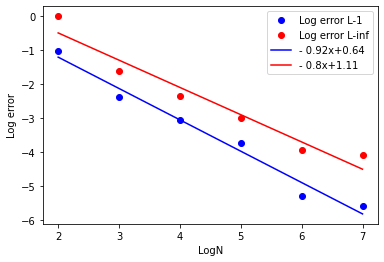}
	\caption{Convergence rate on $N$ for model Example \ref{ex mulitplicative model} with parameters $d=10, \sigma=0.3$ and $K=20$, $J=400$.}
\end{figure}


Finally, in Table \ref{ta comp time}, we report the computational time and the error associated to different dimensions $d = 1,5,10$ in Example \ref{ex mulitplicative model} where $\sigma=0.3$ with NNs $\&$ Upwind scheme ($K=20$, $J=400$) and with fixed number of time step for the splitting $N=32$. Per our specification, the computational time does not increase exponentially and, importantly, neither the empirical error. This behavior is expected from the non-linear regression using neural networks.
\begin{center}
	\begin{tabular}{|c | c| c |c |} 
		\hline
		Dimension & $d=1$& $d=5$ &$d=10$\\ [0.5ex] 
		\hline\hline
		Time & 673s & 1077s & 1488s\\
	\hline
	$L^{1}$ Error & 0.0431  & 0.0594  & 0.0756\\
	\hline
	$L^{\infty}$ Error & 0.0867  & 0.1041 & 0.1253 \\
	\hline
	\end{tabular}
	\captionof{table}{Computational cost in example \ref{ex mulitplicative model} for different dimension $d$ (for the $P$-variable).}
\label{ta comp time}
\end{center}

\bibliographystyle{plain}

\bibliography{biblio}

\end{document}